%%AmsTeX file

%%%Some references updated February 17, 2002
%%%Corrections August 20, 2001
%%Final version, submitted June 6, 2001
%%Corrected misprints May 2002
\documentstyle{amsppt}
%\pagewidth{6in}\vsize8.5in\parindent=6mm\parskip=3pt\baselineskip=14pt
%\tolerance=10000\hbadness=500
%\magnification=1100
%\NoRunningHeads
%\NoBlackBoxes
%\NoLogo
\loadbold
\topmatter
\title
Oscillatory and Fourier Integral operators
 with degenerate canonical relations
\endtitle
\author Allan Greenleaf\qquad\qquad\qquad  Andreas Seeger
\endauthor
\thanks
Supported in part by
grants from the National Science Foundation.
The second author would like to thank the organizers of
the 2000 El Escorial Conference in
Harmonic Analysis and Partial Differential Equations 
for a very interesting
meeting and for  the opportunity
to present some of this material.
\endthanks
\address
University of Rochester, Rochester, NY 14627
\endaddress
\address
University of Wisconsin, Madison, WI 53706
\endaddress
\leftheadtext{A. Greenleaf and A. Seeger}
\rightheadtext{Integral operators with  degenerate canonical relations}
%\date \enddate
\endtopmatter
\document

\define\vectell{\vec l}
\def\tilI{{\widetilde I}}
\define\ic{\imath}
\define\innmod#1#2{#1\!\cdot\!#2}
\define\pdr#1#2{{\frac{\partial^{#2}}{(\partial #1)^{#2}}}}
\define\pder#1#2#3{{\frac{\partial^{#2}{#3}}{(\partial #1)^{#2}}}}

%%x\define\n{{\noindent}}

\define\M{\Cal M}

\define\vth{\vartheta}
\define\R{{\Bbb R}}

\define\cf{{\it cf}}

\define\rank{{\text{\rm rank }}}

\redefine\ker{{\text{\rm Ker }}}

\define\dist{{\text{\rm dist}}}

\define\supp{{\text{\rm supp }}}

\define\inn#1#2{\langle#1,#2\rangle}
\define\rta{\rightarrow}

\define\lcontr{\rfloor}
\define\lco#1#2{{#1}\lcontr{#2}}
\define\lcoi#1#2{\imath({#1}){#2}}
\define\rco#1#2{{#1}\rcontr{#2}}
\redefine\exp{{\text{\rm exp}}}
\define\bin#1#2{{\pmatrix {#1}\\{#2}\endpmatrix}}

\define\ad{{\text{\rm ad}}}

\define\lc{\lesssim}
\define\gc{\gtrsim}

%Greek letters
             
\define\bet{\beta}

\define\eps{\varepsilon}

\define\la{\lambda}             
             
\define\si{\sigma}              
\define\vphi{\varphi}
             
\define\om{\omega}              \define\Om{\Omega}

\define\fA{{\frak A}}

\define\fN{{\frak N}}

\define\fS{{\frak S}}
\define\fT{{\frak T}}

\define\fX{{\frak X}}
\define\fY{{\frak Y}}

\define\fg{{\frak g}}

%Attn:\fi can't be defined.

\define\fm{{\frak m}}

\define\bbM{{\Bbb M}}

\define\bbR{{\Bbb R}}

\define\cA{{\Cal A}}

\define\cC{{\Cal C}}
\define\cD{{\Cal D}}

\define\cF{{\Cal F}}

\define\cH{{\Cal H}}

\define\cK{{\Cal K}}

\define\cM{{\Cal M}}
\define\cN{{\Cal N}}

\define\cP{{\Cal P}}

\define\cR{{\Cal R}}

\define\cT{{\Cal T}}
\define\cU{{\Cal U}}

\define\cW{{\Cal W}}

%roman letters with a tilde

\def\Itil{{\widetilde I}}

%roman letters with a bar
%\define\abar{{\bar a}} \define\Abar{{\bar A}}
%\define\bbar{{\bar b}} \define\Bbar{{\bar B}}
%\define\cbar{{\bar c}} \define\Cbar{{\bar C}}
%\define\dbar{{\bar d}} \define\Dbar{{\bar D}}
%\define\ebar{{\bar e}} \define\Ebar{{\bar E}}
%\define\fbar{{\bar f}} \define\Fbar{{\bar F}}
%\define\gbar{{\bar g}} \define\Gbar{{\bar G}}
%\define\hBar{{\bar h}} \define\Hbar{{\bar H}}
%\define\ibar{{\bar i}} \define\Ibar{{\bar I}}
%\define\jbar{{\bar j}} \define\Jbar{{\bar J}}
%\define\kbar{{\bar k}} \define\Kbar{{\bar K}}
%\define\lbar{{\bar l}} \define\Lbar{{\bar L}}
%\define\mbar{{\bar m}} \define\Mbar{{\bar M}}
%\define\nbar{{\bar n}} \define\Nbar{{\bar N}}
%\define\obar{{\bar o}} \define\Obar{{\bar O}}
%\define\pbar{{\bar p}} \define\Pbar{{\bar P}}
%\define\qbar{{\bar q}} \define\Qbar{{\bar Q}}
%\define\rbar{{\bar r}} \define\Rbar{{\bar R}}
%\define\sbar{{\bar s}} \define\Sbar{{\bar S}}
%\define\tbar{{\bar t}} \define\Tbar{{\bar T}}
%\define\ubar{{\bar u}} \define\Ubar{{\bar U}}
%\define\vbar{{\bar v}} \define\Vbar{{\bar V}}
%\define\wbar{{\bar w}} \define\Wbar{{\bar W}}
%\define\xbar{{\bar x}} \define\Xbar{{\bar X}}
%\define\ybar{{\bar y}} \define\Ybar{{\bar Y}}
%\define\zbar{{\bar z}}  \define\Zbar{{\bar Z}}

%roman letters with a hat

\define\Ga{\Gamma}
\define\lrg{{\text{large}}}
\define\sm{{\text{small}}}

\define\br{\Bbb R}
\define\dg{\dot\gamma}

\define\pr{\pi_R}
\define\pd{\partial}
\define\p#1{\psi^{(#1)}}
\define\g{\gamma}

\define\a{\alpha}
%\define\b{\beta}
%\define\bin#1#2{\Bigl(\buildrel{#1}\over{#2}\Bigr)}
\define\hX#1{\widehat{X}_{#1}}
\define\tc{\tilde c}

We shall mostly survey
%give an overview about recent
results concerning the $L^2$ boundedness of
oscillatory and Fourier integral operators.
Many mathematicians have contributed important results to this subject.
This article does not intend to  give a broad overview; it mainly
focusses on a few topics directly related to the work of the authors.

\head{ \bf 1. The nondegenerate situation}\endhead
\subheading{1.1.  Oscillatory integral operators}
The main subject of the article concerns
oscillatory integral operators given by
$$
T_\la f(x)=\int  e^{\ic\lambda\Phi (x,y)} \sigma (x,y) f (y) dy.
\tag 1.1
$$
In (1.1) it is assumed that the real-valued phase function $\Phi$
 is smooth in $\Omega_L\times \Omega_R$ where $\Omega_L, \Omega_R$
are open subsets of $\bbR^d$ and
amplitude $\sigma \in C^\infty_0( \Omega_L \times\Omega_R)$.
(The assumption that $\dim(\Omega_L)=\dim(\Omega_R)$ is only for
convenience; many of the
definitions, techniques and results described below have some analogues in
the
non-equidimensional setting.)

The $L^2$ boundedness   properties of $T_\la$ are determined
by  the geometry of the canonical relation
$$C = \{ ( x,\Phi_x, y, -\Phi_y ): (x,y)\in \supp \sigma \}
\subset T^*\Om_L\times T^*\Om_R.$$
The best possible situation
occurs when $C$ is locally the  graph of a canonical transformation;
i.e., the projections
$\pi_L$, $\pi_R$ to $T^*\Om_L$, $T^*\Om_R$, resp.,

$$\gathered
\\
C \\
\swarrow \quad \quad
%\downarrow
\quad \quad \searrow \\
T^*\Omega_L
\quad\quad\quad\,  \quad \quad\quad
T^*\Omega_R   \\ \\
\endgathered
$$
%$$\matrix
%      &    &  C  &   &\\
%      &  \swarrow &   & \searrow &  \\
%     &&&&\\
%  &  T^* \Omega_L \quad  &     &\quad   T^* \Omega_R   &
%\endmatrix$$
are locally diffeomorphisms.
In this case H\"ormander \cite{37},\cite{38} proved that the  norm of
$T_\la$ as a bounded operator
on $L^2(\bbR^d)$ satisfies
$$
  \| T_\lambda \|_{L^2 \rightarrow L^2} = O (\lambda^{-d/2}).
\tag 1.2
$$
The proof consists in  applying Schur's test to the kernel
of $T_\la^*T_\la$; see the argument following (1.6)
below.

It is also useful to study a more general class of oscillatory integrals
which naturally arises  when composing
two different operators $T_\la$, $\widetilde T_\la$ and which is also
closely related to the concept of  Fourier integral operator.
We consider  the oscillatory integral  kernel
%operator
 with frequency variable $\vartheta\in \Theta$ (an open subset of
$\bbR^N$),
defined by
$$
K_\lambda (x,y) = \int e^{\ic \lambda \Psi (x,y,\vartheta)}
a(x,y,\vartheta )d\vartheta
\tag 1.3
$$
where  $ \Psi \in C^\infty(\Om_L\times\Om_R\times\Theta)$ is real-valued
 and
$ a \in C^\infty_0(\Om_L\times\Om_R\times\Theta)$. Let $\fT_\la$ be the
associated integral operator,
$$\fT_\lambda f(x) = \int K_\lambda (x,y) f(y) dy.
\tag 1.4$$

Again the $L^2$ mapping properties
of $\fT_\la$ are determined by the geometric properties of the canonical
relation
$$
C = \{ (x, \Psi_{x},y, - \Psi_y) : \Psi_\vartheta = 0 \} \subset
T^* \Omega_L \times T^* \Omega_R .
$$
It is always assumed that $C$ is an immersed manifold, which is a
consequence of
the linear independence of the vectors
$\nabla_{(x,y,\vartheta)} \Psi_{\vartheta_i}$, $i=1,\dots,N$ at
$\{\Psi_\vartheta=0\}$.
In other words,
$\Psi$ is a
nondegenerate phase in
the sense of H\"ormander \cite{37}, although $\Psi$ is not
 assumed to be homogeneous.

As before, the best possible situation for $L^2$ estimates arises
when $C$ is locally the graph of a canonical transformation.
Analytically this means that

$$\det \left ( \matrix \Psi_{xy}   & \Psi_{x\vartheta }\\
                      \Psi_{\vartheta y} & \Psi_{\vartheta
\vartheta}\endmatrix \right ) \not = 0
\tag  1.5
$$

Under this assumption the $L^2$ result becomes
$$\|\fT_\la\|_{L^2\to L^2}\lc \la^{-(d+N)/2}
\tag 1.6$$
so that we discover (1.2) when $N=0$. The proof of (1.6)
could be given by using methods in \cite{37}
or alternatively by a straightforward modification of the argument in
\cite{38}. Indeed
consider the Schwartz kernel $H_\la$  of the operator $\fT_\la^*\fT_\la$
which is given by
$$
H_\lambda (u,y) = \iiint e^{-\ic \lambda [ \Psi (x,u,w) - \Psi
(x,y,\vartheta)]} \gamma(x,u,w,y,\vartheta) dw d\vartheta dx
$$
where $\gamma$ is smooth and compactly supported. By using partitions of
unity we may assume that $\sigma$ in (1.1)  has small support;
 thus  $\gamma$ has small support. Change variables
 $w=\vartheta+h$, and, after interchanging the order of integration,
 integrate parts with respect to the  variables $(\vartheta,x)$. Since
$$
\aligned
\nabla_{x,\vth}
 [ \Psi (x,u,\vartheta + h) - \Psi
(x,y,\vartheta )] =
              \left ( \matrix \Psi_{xy}  & \Psi_{x\vartheta } \\
                    \Psi_{\vartheta y}   & \Psi_{yy} \endmatrix \right
              ) \left ( \matrix u - y\\
                             h \endmatrix \right ) + O(|u-y|^2+|h|^2)
\endaligned
$$
this yields, in view of the small support of $\gamma$,
$$\aligned
 |K_\lambda (u,y) | & \lesssim \int (1 + \lambda |u - y|
              + \lambda |h|)^{-2M} dh \\
         & \lesssim \lambda^{-N -d} {\lambda^d \over (1 + \lambda |u-y|)^M}
\endaligned
$$
if $M>d$.
It follows that $
\| \fT_\lambda^* \fT_\lambda \|_{L^2\to L^2}
 \lesssim \lambda^{-N-d}$ and hence (1.6).

\subheading{1.2 Reduction of frequency variables}

Alternatively, as in the theory of Fourier integral operators, one may
compose $T_\la$
with  unitary operators associated to canonical transformations, and
together  with stationary phase calculations, deduce estimates for operators
of the form (1.3-4) from operators of the form (1.1), which involve no
 frequency variables; in fact this
procedure turns out to
be very useful when estimating operators with
degenerate canonical relations.

We briefly describe the idea based on \cite{37}, for details see
\cite{25}.

Consider the operator  $\fT_\lambda$ with kernel
%be the operator with kernel
$\int_{\bbR^N} e^{\ic \lambda \phi (x,y,z) } a(x,y,z) dz.$  Let
$A_i$,  $i=1,2$,  be symmetric $d\times d$ matrices   and define
$$S_{\lambda ,i} g(x) = \left ( {\lambda \over 2\pi} \right )^{d/2} \int
e^{-\ic  \lambda [ \langle x,w\rangle + {1\over 2} A_i w \cdot w] } g(w) dw;
$$
clearly
$S_{\lambda,i}$ are unitary operators on $L^2(\Bbb R^d)$.
A computation yields that  the operator
$\la^{-d}S_{\lambda ,1} \fT_\lambda S^*_{\lambda ,2}$ can be written as
the sum of an oscillatory integral operator with kernel $O_\la(x,y)$
plus an operator with $L^2$ norm $O(\la^{-M})$ for any $M$.
The oscillatory kernel $O_\la(x,y) $ is again of  the form (1.3) where the
phase function
is given by
$$\Psi(x,y,\vartheta)=
\inn{y}{\tilde w}-\inn {x}{w}+
\frac 12(A_1\tilde w\cdot\tilde w-A_2w\cdot w)+\phi(w,\tilde w,z)
$$ with frequency variables   $\vartheta = (w,z,\tilde w)
\in\Bbb R^d \times \Bbb R^N
\times \Bbb R^d$, and the amplitude is
compactly supported.

  One can choose $A_1,
A_2$ so that for tangent vectors  $\delta x,\delta y \in
\Bbb R^d$  at a reference point the
vector $(\delta x, A_1 \delta x, \delta y, A_2 \delta y)$ is
tangent to the canonical relation $\widetilde C$
associated with $S_{\lambda ,1}
\fT_\lambda S_{\lambda ,2}^*$. Let $\pi_{\text{space}}$ be the projection
$\widetilde C\to \Omega_L\times \Omega_R$ which with our choice of $A_1$,
$A_2$ has invertible differential.
Since the number of frequency variables $(N+2d)$ minus the rank of
$\phi_{\vth \vth}$ is equal to $2d- \rank d\pi_{\text{space}}$, we deduce
that
$\det\phi_{\vth\vth}\neq 0$.

In the integral defining the kernel of
$S_{\lambda ,1}
\fT_\lambda S_{\lambda ,2}^*$
we can now apply the method of stationary phase to reduce the
number of frequency variables to zero, and gain a factor of
$\lambda^{-(2d+N)/2}$.
Thus we may write
$$
S_{\lambda ,1}  \fT_\lambda S_{\lambda ,2}^* = \lambda^{- {N/ 2}}
T_\lambda + R_\lambda
$$
where $T_\la$ is an oscillatory integral operator (without frequency
variables)
and $R_\lambda$ is an operator with $L^2$ norm $O( \lambda^{-M} ) $ for any
large $M$. Since $S_{\la,i}$ are unitary the $L^2$ bounds  for
$\la^{N/2}\fT_\la$ and $T_\la$ are equivalent.

\medskip

\subheading{1.3 Fourier integral operators}
The kernel of a Fourier integral operators  $\cF:
C^\infty_0(\Omega_R)\to \cD'(\Omega_L)$
of order $\mu$, $\cF\in I^{\mu}(\Omega_L, \Omega_R;\cC)$
is locally given as a finite sum of oscillatory integrals
$$
\int e^{\imath \Psi(x,y,\theta)} a(x,y,\theta) d\theta,
\tag 1.7
$$
where now $\Psi$ is nondegenerate in the sense of H\"ormander \cite{37},
satisfies
the homogeneity condition $\Psi(x,y, t\theta)=t\Psi(x,y,\theta)$ for
$|\theta|=1$ and $t\gg 1$, and  $a$ is a symbol of order $\mu+(d-N)/2$.
We assume in what follows that
$a(x,y,\theta)$ vanishes for $(x,y)$ outside a fixed compact set.
The canonical relation is locally given by
$\cC=\{(x,\Psi_x,y,-\Psi_y), \Psi_\theta=0\}$ and we assume that
$$\cC\subset \Bigl(T^*\Om_L\setminus 0_L\Bigr)\times
\Bigl(T^*\Om_R\setminus 0_R\Bigr),
$$
where $0_L$, $0_R$ denote the zero-sections in
$T^*\Om_L$ and $T^*\Om_R$. Staying away from the zero sections implies
$$|\Psi_x(x,y,\theta)|\approx |\theta|\approx|\Psi_y(x,y,\theta)|
\tag 1.8$$
for large $\theta$ (when $\Psi_\theta$ is small).
Let $\beta\in C^\infty_0(1/2,2)$ and
$$a_k(x,y,\theta)=\beta(2^{-k}|\theta|)a(x,y,\theta)$$ and let
$\cF_k$ be the dyadic localization of $\cF$; i.e. (1.7) but with $a$
replaced by
$a_k$. The assumptions $\Psi_x\neq 0$ and $\Psi_y\neq 0$
can be used to show that for $k,l\ge 1$ the operators $\cF_k$ are almost
orthogonal, in the sense that $\cF_k^*\cF_l$ and $\cF_k \cF_l^*$ have
operator norms $O(\min \{2^{-k M}, 2^{-lM}\})$ for any $M$,
provided that $|k-l|\ge C$ for some large but fixed constant $C$.
This follows from  a straightforward integration by
parts
argument based on (1.8) and the assumption of compact $(x,y)$ support.
Using a change of variable $\theta=\la \vth$ the study of the $L^2$
boundedness (and $L^2$-Sobolev boundedness) properties is reduced to
the study of oscillatory integral operators (1.3-4) and, in the
nondegenerate case,
 an application of
estimate (1.2) above.
The result is that if $\cF$ is of order
$\mu$
%%(which means the symbol $a$ in (1.7) is of order $\mu+(d-N)/2$)
and if the associated homogeneous canonical transformation  is  a local
canonical graph,
then $\cF$ maps the Sobolev space $L^2_\alpha$ to $L^2_{\alpha-\mu}$

An important subclass is the class of conormal operators
 associated to  phase functions linear in the frequency variables
(see \cite{37,\S2.4}).  The  generalized Radon transforms
$$
\cR f(x)= \int_{\cM_x} f(y)\chi(x,y)d\sigma_x(y)
\tag 1.9
$$
arise as model cases. Here
$\cM_x$ are codimension $\ell$ submanifolds in $\Bbb R^d$, and $d\sigma_x$
is a smooth
density on $\cM_x$, varying smoothly in $x$, and
$\chi\in C^\infty_0(\Om_L\times\Om_R)$.
One assumes that the $\cM_x$ are sections of a manifold
$\cM\subset\Om_L\times\Om_R$,
 so that the projections to $\Omega_L$ and to
$\Omega_R$ have surjective differential; this assumption insures
the $L^1$ and $L^\infty$ boundedness of the operator $\cR$.
We refer to $\cM$ as the associated
 {\it incidence relation}.

Assuming that $\cM$ is given by an $\bbR^\ell$ valued defining function
$\Phi$,
$$\cM=\{(x,y):\Phi(x,y)=0\},\tag 1.10
$$ then
the distribution kernel of $\cR$ is $\chi_0(x,y)\delta(\Phi(x,y))$
where
$\chi_0\in C^\infty_0(\Om_L\times\Om_R)$ and $\delta$ is the Dirac measure
in $\Bbb R^\ell$ at the origin.
The assumptions on the projections to $\Omega_L$, $\Omega_R$ imply
 that $\rank \Phi_x=\rank \Phi_y =\ell$ in a neighborhood of
$\cM=\{\Phi=0\}$.
The Fourier integral description is then
obtained by writing out $\delta$ by means of the Fourier inversion formula
in $\bbR^\ell$,
$$\chi_0(x,y)\delta(\Phi(x,y))=\chi_0(x,y) (2\pi)^{-\ell}\int_{\Bbb R^\ell}
e^{\ic\tau\cdot\Phi(x,y)} d\tau;
\tag 1.11
$$
this has been used in \cite{35}
where $\cR$
is identified as a Fourier integral operator of order $-(d-\ell)/2$, see also
 \cite{55}.
 More general conormal operators are obtained
by composing Radon transforms
with pseudo-differential operators (see \cite{37}).

The canonical relation associated to the generalized Radon transform
 is the twisted conormal
 bundle of the incidence relation,
$$\Cal C=N^*\cM'=\{(x,\tau\cdot\Phi_x, y,-\tau\cdot\Phi_y):
\Phi(x,y)=0\}.\tag 1.12$$
We can locally (after possibly a change of coordinates) parametrize $\cM$ as
a graph so that
$$\Phi(x,y)=S(x,y')-y''\tag 1.13 $$
with $y'=(y_1,\dots,y_{d-\ell})\in \bbR^{d-\ell}$,
$y''=(y_{d-\ell+1},\dots,y_d)\in \Bbb R^{\ell}$, $S=(S^1,\dots,S^\ell)$.
Using (1.5) with  $\Psi(x,y,\tau)=\tau\cdot\Phi(x,y)$
one verifies that the condition for
$N^*\cM'$ being a local canonical graph is equivalent to the nonvanishing of
the determinant
$$
\det\pmatrix
\tau\!\cdot\!\Phi_{xy}&\Phi_x
\\
{}^t\Phi_y&0\endpmatrix = (-1)^\ell
\det\pmatrix
\innmod{\tau}{S_{x'y'}} &S_{x'}
\\
\innmod{\tau}{S_{x''y'}} &S_{x''}
\endpmatrix
\tag 1.14
$$
for all $\tau \in S^{\ell-1}$.
Under this  condition $\cR$ maps $L^2$ to $L^2_{(d-\ell)/2}$.

We note that the determinant in (1.14) vanishes
for some $\tau$ if $\ell<d/2$. In particular if $\ell=d-1$ then the
expression (1.14) is a linear functional of $\tau$
and thus, if $(x,y)$ is fixed, it
  vanishes for all $\tau$ in a hyperplane. Therefore degeneracies always
occur for averaging over manifolds with high codimension, in particular for
curves in three or more dimensions.

\head{\bf 2. Finite type conditions}\endhead

\subheading{ 2.1. Finite type}
Different  notions of finite type are useful in different situations. Here
we shall restrict ourselves to maps (or pairs of maps) which have
corank $\le 1$.

Let $M$, $N$ be  $n$-dimensional manifolds, $P\in M $ and
$Q\in N $, and let  $f:M \to N $ be a $C^\infty$ map with $f(P_0)=Q_0$.
A vector field $V$ is a {\it kernel  field}
 for the map $f$ on a neighborhood
$\cU$ of $P_0$  if $V$
is smooth on $\cU$
 and
if   there exists a smooth vector field $W$ on $f(U)$ so that
$Df_P V= \det (Df_P)W_{f(P)}$ for all $P\in \cU$. If
  $\rank Df_{P_0}\ge n-1$ then it is easy to see that
 there  is a neighborhood of
$P$ and a nonvanishing kernel vector field $V$ for $f$ on $\cU$. Moreover if
$\widetilde V$ is another kernel field on $U$
 then $\widetilde V= \alpha V-\det(Df)W$
in some neighborhood of $P_0$, for some vector field $W$ and smooth function
$\alpha$.
If $Df=\pmatrix A&b\\c^t&d\endpmatrix$ with $A$ an invertible  $(n-1)\times
(n-1)$ matrix, then $\det Df= \det A(d-c^tA^{-1}b)$ and a choice for the
kernel vector field is
$$V=\frac{\partial}{\partial x_d}- A^{-1}b\cdot \nabla_{x'}.
\tag 2.1$$

\definition{Definition}  Suppose that
$M $ and $N$ are smooth $n$-dimensional manifolds and
that
 $f:M \to N $ is a smooth map with   $\dim \text{ker}( Df) \le 1$
on $M $.
We say that $f$ is {\it of type $k$ at $P$} if
there is a nonvanishing kernel field $V$ near $P$
 so that $V^j(\det Df)_P=0$ for $j<k$ but
$V^k(\det Df)_P\neq 0$.
\enddefinition

This definition was proposed by  Comech \cite{13}, \cite{15}
who assumes in addition that $Df$ drops rank  simply
 on the singular variety $\{ \det Df=0\}$.

The finite type condition  is satisfied for the class of Morin singularities
(folds, cusps, swallowtails, ...) which we shall now discuss.

\subheading{2.2 Morin singularities }
We consider as above maps $f:M\to N$ of corank $\le 1$.
We say that
 $f$  drops rank simply at $P_0$ if
$\rank Df_{P_0}= n-1$ and if
$d(\det Df)_P\neq 0$. Then near $P_0$ the variety
$S_1(f)=\{x: \rank Df =n-1\}$ is  a hypersurface and we say that
 $f$ has an
$S_1$
singularity at $P$ with
singularity manifold $S_1(f)$.

Next let $\fS$ be a hypersurface in a manifold $\cU$ and let
$V$ be a vector field defined on $\fS$ with values in
$T\cU$ (meaning that $v_P\in T_P\cU$ for $P\in \fS$).
We say that $v$ is {\it transversal} to $\fS$ at $P\in\fS$ if $v_P\notin
T_P\fS$.
 We say that $v$ is
{\it simply tangent} to $\fS$ at $P_0$ if there is a one-form $\omega$
annihilating  vectors tangent to $\fS$ so that
$\inn{\om}{v}\big|_\fS$  vanishes of exactly first order at $P_0$.
This  condition does not depend on the particular choice of $\om$.
 Next let
 $P\to \ell(P)\subset T_P(V)$ be a smooth field of lines defined on
$\fS$. Let $v$ be a nonvanishing
 vector field so that $\ell(P)=\Bbb R v_P$.
The definitions of transversality and simple tangency carry over
to field of lines (and the notions
 do not depend on the particular choice of the vector field).

Next consider $F:\cU\to N$ where $\dim \cU=k\ge 2$ and $\dim N=n\ge k$
and assume that
$\rank DF\ge k-1$.
Suppose that $\fS $ is a hypersurface in $\cU$
such that $\rank DF= k-1$ on $\fS$. Suppose that $\ker DF$ is  simply
tangent  to $\fS$ at $P\in\fS$.
Then there is a neighborhood $U$ of $P$ in
$\fS$ such that  the variety
$\{Q\in U: \rank DF\Big|_{T_Q\fS}=k-2\}$
is a smooth hypersurface in $\fS$.

With these notions we can now recall the definition of
Morin singularities (\cite{78}, \cite{47}).
\definition{ Definition } Let $1\le r\le n$.
Let $\fS_1,\dots,\fS_r$ be  submanifolds of an open set $\cU\subset M$
so that
$\fS_k$ is  of dimension $n-k$ in $V$  and
$\fS_1\supset\fS_{2}\supset\dots\supset\fS_r$; we also set $\fS_0:= \cU$.

 We say that $f$ has an
 $S_{1_r}$ singularity  in $\cU$, with a
descending flag of
 singularity
manifolds $(\fS_1,\dots,\fS_r)$ if the following conditions hold in $\cU$.

\roster

\item"{(i)}"
 For $P\in \cU$, either $Df_P$ is bijective or $f$ drops rank simply at $P$.

\item"{(ii)}" For
$1\le i\le r$,    $\rank D(f\big|_{\fS_{i-1}})_P=n-i+1$ for
all $P\in \fS_{i-1}\setminus\fS_i$.

\item"{(iii)}"
 For $2\le i\le r-1$, $\ker   D(f\big|_{\fS_{i-1}})$ is
simply tangent to $\fS_i$ at points in
$\fS_{i+1}$.
\endroster
\enddefinition

\definition{Definition}
We say that $f$ has an  $S_{1_r,0}$ singularity at $P$,
if the following conditions hold.

\roster
\item"{(i)}" There exists a neighborhood
$\cU $ of $P$
 submanifolds $\fS_k$  of dimension $n-k$ in $U$  so that
$P\in \fS_r\subset\fS_{r-1}\subset\dots\subset\fS_1$ and so that
$f:\cU\to N$
 has an
 $S_{1_r}$ singularity in $\cU$,
 with singularity
manifolds $(\fS_1,\dots,\fS_r)$.

\item"{(ii)}"
  $\ker Df_P\cap
T_P(\fS_r)=\{0\}$.
\endroster
\enddefinition

The singularity
manifolds $\fS_k$ are denoted by
 $S_{1_k}(f)$ in singularity theory
(if the neighborhood is understood).
An $S_{1,0}$ (or $S_{1_1,0}$) singularity  is a Whitney fold; an
$S_{1,1,0}$ (or $S_{1_2,0}$) singularity is referred
 to as a Whitney or simple cusp.

If $f$ is given in {\it adapted coordinates} vanishing at $P$, i.e.
$$f: t\mapsto (t', h(t))\tag 2.2$$
then $f$ has an $S_{1_r}$ singularity
in a neighborhood of $P =0$ if and only if
$$(\partial/\partial t_n)^k h(0)=0, \qquad 1\le k\le r,
\tag 2.3$$
and
the gradients
$$\nabla_t\big(\frac{\partial^k h}{\partial t_n^k}\big), k=1,\dots r-1,
\tag 2.4$$ are linearly
independent  at $0$.
Moreover $f$ has an $S_{1_r,0}$ singularity  at $P$ if in addition
$$(\partial/\partial t_n)^{r+1}h(0)\neq 0.
\tag 2.5$$
The singularity manifolds are then given by
$$S_{1_k}(f)=\{t:
(\partial /\partial t_n)^j f(t)=0, \, 1\le j\le k
\}.$$
In these coordinates the kernel field for $f$ is
$\partial/\partial t_n$ and the map  $f$ is of type $r$
at $P$.

Normal forms of $S_{1_r}$ singularities are due to Morin \cite{47}, who
showed that there exists adapted  coordinate systems so that
(2.2) holds with
$$h(t) = t_1t_n+t_2t_n^2+\dots+t_{r-1}t_n^{r-1}+t_n^{r+1}.
\tag 2.6
$$

Finally we mention the situation of maximal degeneracy for
$S_1$ singularities which occurs when the kernel of $Df$ is
everywhere  tangential to the singularity surface $S_1(f)$. In this case
we say that $f$ is a {\it blowdown}; see example 2.3.3 below.

\subheading{2.3. Examples}  We now discuss some
 model examples.
The first set of examples concern translation invariant averages over
curves, the second set
restricted X-ray transforms for rigid line complexes. The map $f$ above will
always be one of the projections $\pi_L:\cC\to T^*\Om_L$ or $\pi_R:\cC\to
T^*\Om_R$. Note that
$S_1(\pi_L)=S_1(\pi_R)$.

{\bf 2.3.1.} Consider the operator on functions in $\Bbb R^d$
$$\cA f(x)=\int f(x+\Gamma(\alpha)) \chi(\alpha)d\alpha
\tag 2.7$$
where $\alpha\to \Gamma(\alpha)$ is a  curve in $\Bbb R^d$ so that
$\Ga'(\alpha), \Gamma''(\alpha),\dots, \Gamma^{(d)}(\alpha)$ are linearly
independent.
Then the canonical relation is given
by
$$C=\{(x, \xi; x+\Gamma(\alpha), \xi): \inn{\xi}{\Gamma'(\alpha)}=0\}.$$
Consider the projection $\pi_L$ then it is not hard to see that
$S_{1_k}(\pi_L)$ is the submanifold of $\cC$ where
in addition
$\inn{\xi}{\Gamma^{(j)}(\alpha)}=0\}$ for $2\le j\le k+1$. Clearly
then $S_{1_{d-1}}(\pi_L)=\emptyset$ so that we have
an
$S_{1_{d-2,0}}$ singularity. The behavior of $\pi_R$ is of course exactly
the same; moreover
for small  perturbations the projections $\pi_L$ and $\pi_R$ still  have at
most
$S_{1_{d-2,0}}$ singularities. Note that in the translation invariant
setting
we have $S_{1_k}(\pi_L)=S_{1_k}(\pi_R)$, but for
for small variable   perturbations the manifolds
$S_{1_k}(\pi_L)$, $S_{1_k}(\pi_R)$ are typically different if $k\ge 2$.

By Fourier transform arguments
%In either case the $L^2$ Sobolev
 and van der Corput's lemma it is easy to see that
$\cA$ maps $L^2 (\Bbb R^d)$ to the Sobolev-space
$L^2_{1/d}(\Bbb R^d)$ and it is conjectured that this estimate remains true
for variable coefficient perturbations.
This is known in dimensions $d\le 4$ (cf. \S5 below).

{\bf 2.3.2.} Consider the example (2.7) with $d=3$ and $$\Gamma(\alpha)=
(\alpha, \frac{\alpha^m}{m}, \frac{\alpha^n}{n})$$ where $m$, $n$ are
integers with $1<m<n$.

The canonical relation $C$ is given as the set of
$(x,\xi,y,\xi)$ where
$x_2-y_2-(x_1-y_1)^m/m=0$,
$x_3-y_3-(x_1-y_1)^n/n=0$, and $\xi= (\xi_1(\la,\mu),\la, \mu)$ so that
$$\xi_1=-(x_1-y_1)^{m-1}\la-(x_1-y_1)^{n-1}\mu$$
with $(\la,\mu)\neq (0,0)$.

$C$ is thus parametrized by
$(x_1,x_2,x_3,\la,\mu,y_1)$ and the singular variety $S_1(\pi_L)$
is given by the equation
$$(m-1)(x_1-y_1)^{m-2}\la+(n-1)(x_1-y_1)^{n-2}\mu=0.$$
Note that $\partial /\partial y_1$ is a kernel vector field
and hence $\pi_L$ is of type at most $n-2$ everywhere.
Note that $S_1(\pi_L)$ is a  smooth submanifold only if $m=2$. The case
$m=2$, $n=3$ corresponds to the situation considered above
(now $\pi_L$ is a fold).
If $m=2$, $n=4$  we have a simple cusp ($S_{1,1,0}$) singularity and
$S_{1,1}(\pi_L)$ is the submanifold of $S_1(\pi_L)$
on which $x_1=y_1$.
If $m\ge 3$, $n>m$ then the singular variety is not a smooth manifold
but  the union of the
two transverse  hypersurfaces $\{
(m-1)\la+(n-1)(x_1-y_1)^{n-m}\mu=0\}$ and $\{x_1=y_1\}.$

{\bf 2.3.3.} For an example for a one-sided behavior we consider the
restricted X-ray transform
$$
\cR f(x',x_d)=\chi_0(x_d)\int f(x'+t\gamma(x_d),t)\chi(t) dt
\tag 2.8$$
where $\gamma$ is now the regular parametrization of a curve in $\Bbb
R^{d-1}$ and $\chi_0,\chi$ are smooth and compactly supported.
We say that $\cR$ is associated to a $d$ dimensional  line complex
which is referred to as {\it rigid}
because of the translation invariance in the
$x'$ variables.

The canonical relation is now given by
$$C=\big\{
\big(x',x_d, \tau,y_d\tau\cdot\gamma'(x_d);
x'+y_d\gamma(x_d), y_d, \tau,\tau\cdot\gamma(x_d)\big)\big\}
$$
and the singular set $S_1(\pi_L)=S_1(\pi_R)$ is the submanifold on which
$\tau\cdot \gamma'(x_d)=0$. One computes that $V_L=\partial/\partial y_d$ is
a kernel vector field for
$\pi_L$ and
$V_R=\partial/\partial x_d$ is a kernel vector field for $\pi_R$ .
Clearly $V_L$ is tangential to $S_1(\pi_L)$ everywhere so that
$\pi_L$ is a blowdown.
The behavior of the projection $\pi_R$ depends on assumptions on $\gamma$.
The best case occurs when $\gamma'(x_d),\dots, \gamma^{(d-1)}(x_d)$  are
linearly independent everywhere. The singularity manifolds
$\fS_k=S_{1_k}(\pi_R)$ are then given by the equations
$$\tau\cdot\gamma^{(j)}(x_d)=0,
\qquad j=1,\dots, k,$$ and
thus  $S_{1_{d-1}}(\pi_R)=\emptyset$ and
$\pi_R$ has (at most) $S_{1_{d-2},0}$ singularities.

For the model case given here it is easy to derive the sharp $L^2$-Sobolev
estimates.  Observe that
$$R^*R f(w) =\overline{\chi(w_d)} \iint f(w'+s\gamma(\alpha), w_d+s)
|\chi_0(\alpha)|^2 \chi(s) ds\,d\alpha
$$
defines (modulo the cutoff function) a translation invariant operator.
By van der Corput's  Lemma it is easy to see  that
$$\Big|
\iint e^{-\ic s(\xi'\cdot \gamma(\alpha)+\xi_d)}|\chi_0(\alpha)|^2 \chi(s)
ds\,d\alpha
\Big|\lc (1+|\xi|)^{-\frac 1{d-1}}
$$
and one deduces that $\cR$ maps $L^2$ to $L^2_{1/(2d-2)}$.

It is conjectured that the $X$-ray transform for general well-curved line
complexes

$$\cR f(x',\alpha)=\chi(x',\alpha)\int f(x'+s\gamma(x',\alpha),s)\chi(s) ds
\tag 2.9$$
satisfies locally  the same estimate; here the support of $\chi$ is
supported in
$(-\eps,\eps)$ for small $\eps$ and it is assumed that for each fixed $x'$
the vectors $(\partial/\partial\alpha)^j\gamma$, $j=1,\dots, d-1$
 are linearly independent.
The sharp
$L^2\to L^2_{1/(2d-2)}$ estimate is currently known in dimension $d\le 5$
({\it cf.} \S4-5 below).

\subheading{2.4 Strong Morin singularities}

We now discuss the  notion of  {\it strong}  Morin singularities, or
$S_{1_r}^+$
singularities for maps
into  a fiber
bundle
% $\buildrel{\pi_B}\over{W\rta B}$; here
 $W$ over a base manifold $B$, with projection $\pi_B$.
Here
it is assumed that $\dim W=n$ and
$\dim(B)=q\le n-r$, so that the fibers $W_b=\pi_B^{-1}b$ are $n-q$
dimensional manifolds
(see \cite{26}).
The relevant  $W$ is $T^*\Om_R$, the
cotangent bundle of the base $B=\Om_R$.
\definition{Definition } Let $b=\pi_B(f(P))$ and
let $W_{b}=\pi_B^{-1}b$ be the fiber
 through $f(P)$. The map
$f$ has an $S_{1_r,0}^+$ singularity at $P$
 if

(i)  $f$ intersects $W_b$ transversally, so that there is a neighborhood $U$
of $P$ such that the
preimages $f^{-1}W_b\cap U$ are smooth manifolds
of dimension $n-q$,

 and if

(ii)
 $f\big|_{f^{-1}(W_b)\cap U}$  has an  $S_{1_r,0}$ singularity at $P$.
\enddefinition

Now let $\cC\subset T^*\Om_L\times T^*\Om_R$ be a canonical relation,
consider $\pi_L:\cC\to T^*\Om_L$
and use the natural fibration
$\pi_{\Omega_L}:T^*\Om_L\to \Omega_L$.
If  $\pi_L:\cC\to T^*\Omega_L$   has an $S_{1_r,0}^+$
singularity at $c\in C$, $c=(x_0,\xi_0,y_0,\eta_0)$ then near $c$ we can
restrict
$\pi_L$ to $\pi_{\Omega_L}^{-1}(\{y_0\})$
 and define $\pi_{L,y_0}$ as the restriction of $\pi_L$ to
$\pi_{\Omega_L}^{-1}(\{y_0\})$ and $\pi_{L,y_0}$ has an $S_{1_r,0}$
singularity at $c$.

We remark that
for the examples in 2.3.1 both $\pi_L$ and $\pi_R$  have  strong Morin
singularities
while for the example in 2.3.3 $\pi_R$ has strong Morin singularities. This
remains true for small perturbations
of these examples.

In order to verify the occurence of strong Morin singularities  for
canonical relations which
come up in studying averages on curves the following simple lemma is
useful.

\proclaim{Lemma }  Let $I$ be an open interval, let
 $\psi:I\rta\R^n$ be a smooth parametrization of a regular curve not
passing
through $0$ and let
$$M=\{(t,\eta)\in I\times\R^n:\eta\cdot\psi(t)=0,\hbox{ some } t\in I\}.$$
Let $\pi:M\to \Bbb R^n$ be defined by $\pi(t,\eta)= \eta$.

Then
 $\pi$ has
singularities at most $S_{1_{n-2},0}$ if and only if
 $\{\psi(t),\dot\psi(t),\dots,\psi^{(n-1)}(t)\}$
is a  linearly independent set for all $t\in I$.
\endproclaim

For the proof assume first the linear independence of $\psi^{(j)}(t)$. We
may  work near $t=0$ and by a
linear change of variables, we may assume that $\psi^{(j)}(t_0)=e_{j+1},
0\le
j\le
n-1$, where $\{e_j\}_{j=1}^n$ is the standard basis of $\R^n$.
Thus
$$\align\eta\cdot\psi(t)=&\sum_{j=0}^{n-1} \eta_{j+1}
\frac{t^j}{j!}(1+O(|t|))
\cr
=& \eta_1(1+O(|t|))+\sum_{j=2}^{n}
\eta_j\frac{t^{j-1}}{(j-1)!}(1+O(|t|))
\endalign
$$
with $\eta=(\eta_1,\eta')$. We can solve
$\eta\cdot\psi(t)=0$ for $\eta_1=\eta_1(\eta',t)$,
$$\eta_1=-\sum_{j=2}^{n}\eta_j \frac{t^{j-1}}{(j-1)!}(1+O(|t|)).$$
Hence, $(\eta',t)$ and $(\xi',\xi_1)$ form adapted coordinates
(\cf. (2.2))  for the map $\pi$, and in these coordinates
$$\pi(\eta',t)=(\eta',\phi(\eta',t))=
\bigl(\eta',-\sum_{j=2}^{n}\eta_j\frac{t^{j-1}}{(j-1)!}\bigl(1+O(|t|)\bigr)$$
where $\phi$ satisfies
$$\frac{\pd^j\phi}{\pd t^j}(0,0)=0,\quad 1\le j\le n-2,\quad
\frac{\pd^{n-1}\phi}{\pd t^{n-1}}(0,0)\ne0$$
and the differentials
$$\Bigl\{ d(\frac{\pd^j\phi}{\pd t^j})(0,0)\Bigr\}_{j=1}^{n-1}=
\bigl\{e_j\bigr\}_{j=2}^n
$$
are linearly independent.
Thus $\pi$ has at most $S_{1_{n-2},0}$ singularities.

Conversely, assume that
$\pi$ has at most $S_{1_{n-2},0}$ singularities. Since $\psi$ does not pass
through
the origin, we may assume that $\psi_n(t)\neq 0$ locally. Then the map
$\pi$ is given in adapted coordinates by
$$
(\eta',t)\mapsto (\eta', -\sum_{j=1}^{n-1}
\eta_j\frac{\psi_j(t)}{\psi_n(t)})
$$
and the linear independence  follows easily from (2.3-5).

\subheading {2.5 Mixed finite type conditions}
We briefly discuss mixed conditions for pairs of maps $(f_L,f_R) $
where $f_L:M\to N_L$, $f_R:M\to N_R$ where $M$, $N_L$, $N_R$ are all $d$
dimensional and
$f_L$, $f_R$ are {\it volume equivalent}, i.e., there is a nonvanishing
function $\alpha$
so that
$\det Df_L=\alpha \det Df_R$ in the domain under consideration.

Let $V_L$, $V_R$ be nonvanishing kernel fields on $M$ for the maps
$f_L$, $f_R$.
Let $U$ be a neighborhood of $P$ in $M$. We define
$\cD_{j,k}(U)$ to be the linear space of differential operators generated
spanned by operators
of the form
$$a_1V_1\dots a_{j+k} V_{j+k}$$
where $V_i$ are kernel fields for the maps $f_L$ or $f_R$ in $U$,
and $k$ of them are kernel fields for $f_L$ and $j$ of them are kernel
fields for $f_R$.
Let $h$ be a real valued function defined in a neighborhood of $P\in M$; we
say that
$h$ {\it vanishes of order $(j,k)$} at $P$ if $Lh_P=0$ for all
$L\in \cD^{j-1,k}\cup \cD^{j,k-1}$.
We say that $(f_L, f_R)$ is of type $(j,k) $ if
$h\equiv \det Df_L$
vanishes of order $(j,k)$ at $P\in M$
and if there is an operator
$L\in \cD^{j,k}$  so that $Lh_P\neq 0$.
Because of the assumption of volume equivalence $\det Df_L$ in this
definition can be
replaced by $\det Df_R$.
In the   canonical example of interest here we have
  $M=\cC\subset T^*\Omega_L\times T^*\Omega_R$, a canonical relation, and
$f_L\equiv \pi_L$,
$f_R\equiv \pi_R$ are the projections to $T^*\Omega_L$ and $T^*\Omega_R$,
respectively.

\head{\bf 3.  Fourier integral operators in two dimensions}\endhead

In this section we examine the regularity of Fourier integral operators in
two dimensions,
in which case one can get the sharp $L^2$ regularity  properties with the
possible exception of
endpoint
estimates. We shall assume that $\Omega_L$, $\Omega_R$ are open subsets of
$\bbR^2$, $\cC\subset( T^*\Omega_L\setminus 0_L)\times ( T^*\Omega_R\setminus
0_R)$
is a {\it homogeneous} canonical relation and $\cF\in
I^{-1/2}(\Om_L,\Om_R,\cC)$, with compactly supported
distribution kernels; we assume that the rank of the projection
$\pi_{\text{space}}:\cC\to \Om_L\times\Om_R$ is $\ge 2$ everywhere.
 The generalized Radon transform (1.11) (with $\ell=1$, $d=2$)
is a model case in which  $\rank(d\pi_{\text{space}})=3$.

In order to formulate the  $L^2$ results we shall
work with the Newton polygon,
as in
 \cite{58} where oscillatory integral
operators
in one dimension are considered. We recall that for a set
$E$ of pairs $(a,b)$ of nonnegative  numbers the Newton polygon
associated to $E$ is the closed convex hull of all
quadrants $Q_{a,b}=\{(x,y):x\ge a,y\ge b\}$ where $(a,b)$  is taken
from $E$. 

\definition{Definition} For $c\in \cC$ let $\cN(c)$ be the Newton polygon
associated to the set
$$E(c)=\{(j+1,k+1): \cC \text{ is of type } (j,k) \text{ at } c\}.
\tag 3.1 $$
Let $(t_c,t_c)$ the point of intersection of the boundary $\partial \cN(c)$
with
the diagonal $\{(a,a)\}$.
\enddefinition

Using the notion  of type $(j,k)$ in \S 2.5 we can now formulate

\proclaim{3.1. Theorem} Let $\Omega_L, \Omega_R\subset\bbR^2$  and $\cC$ as
above %\footnote{Check}
 and let
$\cF\in I^{-1/2}(\Om_L,\Om_R;\cC)$,
with compactly supported distribution kernel.
Let $\alpha=\min_c (2t_c)^{-1} $.

Then the operator $\cF$ maps $L^2$ boundedly to
$L^2_{\alpha-\eps}$  for all $\eps>0$ .
\endproclaim

In the present two-dimensional situation one can reduce matters
to operators  with  phase functions that are linear in the frequency
variables (i.e., the conormal situation). We briefly describe this
reduction.

First, our operator can be written modulo smoothing operators as a finite
sum of operators of the form
$$
\cF f(x)=\int e^{\ic\vphi(x,\xi)} a(x,\xi) \widehat f(\xi) d\xi
\tag 3.2
$$
where $a$ is of order $-1/2$, and  has compact $x$ support.
We may also assume that $a(x,\xi)$ has $\xi$-support in an annulus
$\{\xi:|\xi|\approx \la\}$ for large $\la$.
By scaling we can reduce matters to show that the $L^2$ operator norm
for the oscillatory integral operator $T_\la$ defined by
$$
T_\la g(x)=\int e^{\ic \la\vphi(x,\xi)}\chi(x,\xi) g(\xi) d\xi
$$
is  $O(\la^{-1/2-\alpha})$; here $\chi$ has compact support and
vanishes for $\xi$ near $0$.
We introduce polar coordinates in the last integral,
$\xi=\sigma(\cos y_1,\sin y_1)$
and put
$$S(x,y_1)=\phi(x_1,x_2,\cos y_1, \sin y_1).$$ Then the asserted bound for
   $\|T_\la\|$
is equivalent to the same bound for the $L^2$ norm of $\widetilde T_\la$
defined by
$$
\widetilde T_\la h(x)=\int e^{\ic\la \sigma S(x,y_1)}
\widetilde \chi
(y_1,\sigma) dy_1 d\sigma
$$
for suitable $\widetilde \chi$;
here we have used the homogeneity of $\vphi$. Now we rescale again
and  apply  a Fourier transform in $\sigma$
and see that the  bound  $\|\widetilde T_\la\|=O(\la^{-1/2-\alpha})$
follows from
the $L^2\to L^2_\alpha$ bound for the
conormal Fourier integral operator with
distribution kernel
$$
\int e^{\ic\tau\Phi(x,y)}
b(x,\tau) d\tau
\tag 3.3
$$
where $\Phi(x,y)=S(x,y_1)-y_2$, and
 $b$ is a symbol of order $0$, supported in  $\{|\tau|\approx \la\}$
and compactly supported in $x$.

Thus it suffices to  discuss conormal operators of this form; in
 fact for them  one can prove
almost sharp $L^p\to L^p_\alpha$ estimates.
Before stating these results we shall first
reformulate the mixed finite type assumption from \S 2.5 in the present
situation.

\subheading{3.2. Mixed finite type conditions in the conormal situation}
We now look at  operators with distribution kernels of the form
(3.3). The singular support of such operators is given by
$$\cM=\{(x,y):\Phi(x,y)=0\}$$ and it is assumed that $\Phi_{x}\neq 0$,
$\Phi_y\neq 0$.
The canonical relation is the twisted conormal bundle $N^*\cM'$ as in
(1.12).
In view of the homogeneity the type condition at
$c_0=(x_0,y_0, \xi_0, \eta_0)\in N^*\cM'$ is equivalent with the type
condition at
$(x_0,y_0, r\xi_0, r\eta_0)$ for any $r>0$ and since the fibers in $N^*\cM'$
 are one-dimensional it seems natural to formulate finite type conditions
in terms of  vector fields tangent to $\cM$, and their commutators. We now
describe these conditions
but refer for a more detailed discussion to \cite{67}. Related ideas have
been used in the  study of subelliptic operators
(\cite{36}, \cite{63}), in complex analysis
(\cite{41}, \cite{2}) and, more recently,  in the study
of singular Radon transforms (\cite{11}).

Two types of vector fields  play a special role:
We say that  a vector field
 $V$ on $\cM$ is of type
$(1,0)$ if $V$ is  tangent to $\cM\cap(\Omega_L\times\{0\})$; likewise we
define $V$ to be of type $(0,1)$ if $V$ is tangent to
$\cM\cap(\{0\}\times\Omega_R)$.
The notation is suggested by an analogous situation
in  several complex variables  (\cite{41},  \cite{55}).

Note that at every point $P\in \cM$
the vector fields of type $(1,0)$ and $(0,1)$ span a two-dimensional
subspace
 of the three-dimensional tangent space $T_P\cM$.
Thus we can pick a nonvanishing $1-$form $\omega$ which
annihilates vector fields of type $(1,0)$ and $(0,1)$; we may choose
$\omega=d_x\Phi-d_y\Phi$ and
$X=\Phi_{x_2}\partial_{x_1}-\Phi_{x_1}\partial_{ x_2}$,
$Y=\Phi_{y_2}\partial_{y_1}-\Phi_{y_1}\partial_{ y_2}$
are $(1,0)$ and $(0,1)$ vector fields, respectively.
With this  choice
$$\inn\om{[X,Y]}=-2
\det\pmatrix
\Phi_{xy}&\Phi_x
\\
{}^t\Phi_y&0\endpmatrix \tag 3.4$$
which is (1.14) in the situation $\ell=1$, $\theta=1$ and relates
 $ \inn\om{[X,Y]}$
to $\det d\pi_{L/R}$. Thus $N^*\cM'$ is a local canonical graph iff
$\inn\om{[X,Y]}$ does not vanish. The quantity (3.4) is often referred to as
``rotational curvature''
(\cf. \cite{55}).

Now let $\mu$ and $\nu$ be two positive integers.
%We now define vector fields of type $(\mu,\nu)$,
%for two pairs of positive integers $(\mu,\nu)$. Let $ad(V)W=[V,W]$.
For a neighborhood $U$ of $P$ let
$\cW^{\mu,\nu}(U)$ be the module generated by vector fields
$\ad W_1\ad W_2\dots \ad W_{\mu+\nu-1}(W_{\mu+\nu})$ where $\mu$ of these
vector fields are of type $(1,0)$ and $\nu$ are  of type $(0,1)$. The
finite type condition in (2.4) can be reformulated as follows.
Let $P\in \cM$ and let $c\in N^*\cM'$ with base point $P$. Then
$\cC$ is of type $(j,k)$ at $c$ if there is a neighborhood $U$ of $P$ so
that
for all vector fields $W\in \cW^{j+1,k}(U)\cup\cW^{j,k+1}(U)$ we have
$\inn {\omega}{ W}_P=0$ but there is a vector field
$\widetilde W$
 in $\cW^{j+1,k+1}$
for which
$\inn {\omega}{ \widetilde W}_P\neq 0$.
\footnote{Here we deviate from the terminology in \cite{67}, where
the incidence relation  $\cM$ is said to be of  type
$(j+1,k+1)$ at $P$.}
Now coordinates can be  chosen so that $\Phi(x,y)=-y_2+S(x,y_1)$
and the generalized
Radon transform  is given by $$\cR f(x)=\int \chi(x,y_1, S(x,y_1))
f(y_1, S(x,y_1)) dy_1
\tag 3.5
$$
where
$S_{x_2}\neq 0$ and $\chi\in C^\infty_0(\Om_L\times \Om_R)$. If
$$\Delta(x,y_1)= \det \pmatrix S_{x_1y_1} &S_{x_1}\\ S_{x_2y_1}&
S_{x_2}\endpmatrix$$
then at $P=(x,y_1, S(x,y_1))$ the mixed finite type condition  amounts to
$$
X^{j'} Y^{k'}\Delta(x,y_1)=0 \text{ whenever $j'\le j$ and $ k'<k$ or
$j'<j$ and $k'\le k$ }
\tag 3.6
$$
but
$$X^{j} Y^{k}\Delta(x,y_1)\neq 0 \tag 3.7$$ for
$X=S_{x_2}\partial_{x_1}-S_{x_1}  \partial_{x_2}$ and
$Y=\partial_{y_1}+S_{y_1}\partial_{y_2}$.
For the equivalence of these conditions see \cite{67}.

We now relate the last condition to the finite type condition above.
Notice that $$\cC=\{(x_1,x_2,\tau S_{x_1}, \tau S_{x_2};y_1, S(x,y_1),
-\tau S_{y_1},\tau)\}$$
and  using coordinates $(x_1,x_2, y_1,\tau)$ a kernel
vector field for the projection $\pi_R$ is given by
$V_R=S_{x_2}\partial_{x_1}-S_{x_1}  \partial_{x_2}$;
 this can be identified
with the vector field $X$ on $\cM$.
Moreover a kernel vector field for the projection $\pi_L$ is given by
$V_L=\partial/\partial_{y_1}-S_{x_2}^{-1}S_{x_2y_1}\partial/\partial \tau$
and  for any function of the  form $F(x,y_1)$ we see that
$(V_L-Y)(\tau F)$ equals $F$ multiplied by a $C^\infty$ function.
Thus it is immediate that
$C$ is of type $(j,k)$ at the point $c$ (with coordinates $(x,y_1,\tau)$) if
conditions (3.6), (3.7) are
satisfied, and this is just a condition at the
base-point $P$.

We shall now return to the proof of Theorem 3.1 and formulate an $L^p$
version for the conormal situation.

\proclaim{3.3. Theorem \cite{67}} Let $\Omega_L, \Omega_R\subset\bbR^2$, let
$\cM\subset \Omega_L\times\Omega_R$ so that the projections to $\Omega_L$
and $\Omega_R$ have surjective differential. Suppose that
$\cF\in I^{-1/2}(\Om_L,\Om_R;N^*\cM')$
with compactly supported distribution kernel.

For $c\in N^*\cM'$ denote by $\widetilde \cN(c)$
 the closure of the image of $\cN(c)$ under
the map $(x,y)\mapsto (\frac{x}{x+y},\frac 1{x+y})$; i.e., the convex hull
of the points
$(1,1)$, $(0,0)$ and $(\frac{j+1}{j+k+2}, \frac 1{j+k+2})$ where
$N^*\cM'$ is of type $(j,k)$ at $c$.

Suppose that $(1/p,\alpha)$ belongs to the interior of $\widetilde \cN(c)$,
for every $c$.
 Then
$\cF$ is bounded from $L^p$ to $L^p_\alpha$.
\endproclaim

The $L^2$ estimate of Theorem 3.1 for conormal operators follows as a
special case,
 and for the general situation we use the above reduction.
Theorem 3.3 is sharp up to the open endpoint cases
(\cf. also \S 3.5.1-3 below).

We now sketch the main ingredients of the proof of
Theorem 3.3. We may assume that $S_{x_2}$ is near $1$ and $|S_{x_1}|\ll 1$.
Suppose that $Q=(x^0,y^0)\in \cM$ and suppose that the type $(j',k')$
condition holds for
 some choice of $(j',k')$ with $j'\le j$ and $k'\le k$ at $Q$, and suppose
that
this type assumption is still valid in a neighborhood on the support of the
cutoff function $\chi$ in (3.5)
(otherwise we work with partitions of unity).

Our goal is then to prove that
 $\cF$ maps
$L^p$ to $ L^{p}_\alpha$ for $p=(j+k+2)/(j+1)$ and $\alpha < 1/(j+k+2)$.

%
%This is done by an interpolation argument which is closely related to the
%idea of damping oscillatory integrals  to get improved $L^2$ estimates,
%see \cite{CowMau}, \cite{SoSt}.$

Since we do not attempt to obtain an endpoint result, it is sufficient to
prove
the required estimate for operators with the frequency variable
localized  to $|\tau|\approx\la$ for large $\la$.
We then make an additional dyadic decomposition in terms of the size of
$|\Delta|$ (i.e., the
rotational curvature). Define  a Fourier integral operator 
$\cF_{\la,l_0}$ by  
$$
\cF_{\la, l_0} f(x)=
\int f(y)
\int  e^{\ic\tau(S(x,y_1)-y_2)}
\beta(x,y,\frac{|\tau|}{\la}) \chi( 2^{l_0}|\Delta(x,y_1)|) d\tau \, dy;
$$ then by interpolation arguments our goal will be
 achieved by proving the following crucial estimates:
$$
\|\cF_{\la,l_0}\|_{L^p\to L^p}\le C_\gamma
2^{-l_0\gamma}, \qquad  p=\frac{j+k}{j}, \quad \gamma<\frac 1{j+k},
\tag 3.8
$$
and
$$
\|\cF_{\la,l_0}\|_{L^2\to L^2}\le C_\eps
2^{l_0(\frac 12 +\eps)}\la^{-1/2}.
\tag 3.9
$$
A variant of this  interpolation argument goes back to
investigations on maximal operators in
\cite{18} and \cite{71}, \cite{72}, and (3.9) can be thought of
 a version of an estimate for damped
oscillatory integrals.

The type assumption is only used for the estimate (3.8).
We note that by integration by parts with respect to the frequency variable
the kernel
of $\cF_{\la,l_0}$ is bounded by
$$\la(1+\la|y_2-S(x,y_1)|)^{-N} \widetilde \chi(2^{l_0}\Delta(x,y_1)).
\tag 3.10$$

We can use a well known sublevel set estimate
related to van der Corput's lemma (see \cite{8}) to see that for each fixed
$x$  the set of all
$y_1$ such that $|\Delta(x,y_1)|\le 2^{-l}$ and
$|\partial^{k'}_{y_1}\Delta(x,y_1)|
\approx 2^{-m}$ has Lebesgue measure
bounded by
$ C_\eps 2^{\eps l_0} 2^{(m-l_0)/{k'}}\le
 C_\eps 2^{\eps l_0} 2^{(m-l_0)/{k}}$ if
$m\le l_0$. Moreover, if $\fS(y,x_1)$ is implicitly defined by
$y_2=S(x_1, \fS(y,x_1),y_1)$ then the assumption
$X^{j'}Y^{k'}\Delta\neq 0$ for some $(j',k')$, $j'\le j$, $k'\le k$ implies
that
$\partial_{x_1}^{j'}[\partial_{y_1}^{k'}\Delta(x_1,\fS(y,x_1),y_1)]\neq 0$,
for
some $(j',k')$, $j'\le j$, $k'\le k$.
Thus  for  fixed $y$ the set of all $x_1$ for which
$|\partial_{x_1}^{j'}[\partial_{y_1}^{k'}\Delta(x_1,\fS(y,x_1),y_1)]|\le
2^{-m}$ has Lebesgue measure $\lc 2^{-m/j'}\lc 2^{-m/j}$.
The two sublevel set estimates together with (3.10)
and straightforward applications
of H\"older's inequality yield (3.8), see \cite{67}.

We now turn to the harder $L^2$ estimate (3.9). We sketch the ideas of the
proof (see
\cite{66} and also  \cite{67} for some corrections).

Firstly,  if $2^{l_0}\le \la$ we consider  as above
the oscillatory integral operator
$T_{\la,l_0} $ given by
$$
T_{\la,l_0} g(x) =
\int  e^{\ic\la\tau(S(x,y_1)-y_2)}
\eta(y_1,\tau) \chi( 2^{l_0}|\Delta(x,y_1)|) g(y_1,\tau) d\tau dy_1
\tag 3.11
$$
with compactly supported $\eta$;
it suffices to show that
$$\|T_{\la,l_0}\|_{L^2\to L^2}\le C_\eps 2^{l_0(1+\eps)/2}\la^{-1}.
\tag 3.12 $$
If $|\Delta(x,y_1)|\le \la^{-1}$ we modify our definition by localizing to
this set.
We note that it suffices to estimate the operator
$\chi_{Q'}\cF[\chi_Q f]$ where $Q$ and $Q'$ are squares  of sidelength
$2^{-l_0 \eps/10}$,
 since summing over all relevant pairs of squares will only introduce an
error
$O(2^{4\eps l_0/10})$ in the final estimate.

If we tried to use the  standard $TT^*$ argument we would have to have good
lower bounds
for $S_{y_1}(w,y_1)-S_{y_1}(x,y_1)$ in the situation where
$S(w,y_1)-S(x,y_1) $ is small, but the appropriate lower bounds fail to hold
if the
rotational curvature  is too small.
 Thus it is necessary to work with finer decompositions.
Solve the equation $S(w,y_1)-S(x,y_1) =0$ by $w_2=u(w_1,x,y_1)$
and  expand
$$\multline
S_{y_1}(w,y_1)-S_{y_1}(x,y_1)=\\
S_{y_1}(w_1,u(w_1,x,y_1),y_1)-S_{y_1}(x,y_1)+O(S(w,y_1)-S(x,y_1))
\endmultline
\tag 3.13
$$
and
$$\multline
S_{y_1}(w_1,u(w_1,x,y_1),y_1)-S_{y_1}(x,y_1)=\sum_{j=0}^M\gamma_j(x,y_1)(w_1
-x_1)^{j+1}
+ O(2^{-l_0\eps/M})
\endmultline
\tag 3.14
$$
where  $M\gg 100/\eps$. In particular
$$\gamma_0 (x,y_1)=S_{y_1x_1}(x,y_1)+u_{w_1}(x_1,x,y_1) S_{y_1,x_2}(x,y_1)=
\frac{\Delta(x,y_1)}{S_{x_2}(x,y_1)};
$$
thus $|\gamma_0(x,y_1)|\approx 2^{-l_0}$. For the coefficient
$\gamma_j(x,y_1)$ we have  $$\gamma_j(x,y_1)= S_{x_2}^{-1}
V^j\Delta(x,y_1)+\sum_{k<j}
\alpha_k(x,y_1) V^k(x,y_1)\Delta(x,y_1)$$ where the $\alpha_k$ are smooth
and $V$
is the $(1,0)$ vector field
$\partial_{x_1}- S_{x_1}/ S_{x_2}\partial_{x_2}$.
We introduce an additional
localization in terms the size of $\gamma_j(x,y_1)$.
For $\vectell=(l_0,\dots,l_M)$, with $l_j<l_0$ for $j=1,\dots, M$
define
$$
T_{\la,\vectell} \,g(x) =\iint
e^{\ic\la\tau S(x,y_1)}\eta(x,y_1)
\prod_{j=0}^M\chi(2^{l_j}|\gamma_j(x,y_1)|) g(y_1,\tau ) dy_1d\tau
$$
which describes a localization to the sets where $|\gamma_j(x,y_1)|\approx
2^{-l_j}$.
A modification of the definition is required if   $|\gamma_j|\le 2^{-l_0}$
for some
$j\in\{1,...,M\}$.

Since we  consider at most  $O((1+l_0)^M)=O(2^{\eps l_0})$ such
operators it suffices to bound any individual $T_{\la,\vectell}$\,, and the
main estimate is

\proclaim{3.4. Proposition}
$$\|T_{\la,\vectell}\,\|_{L^2\to L^2}\le C_\eps
2^{l_0(1+\eps)/2}\la^{-1/2}.$$
\endproclaim

In what follows we fix $\la$ and $\vectell$ and set
$$
\cT=T_{\la, \vectell}\,.$$ The proof of the  asserted $L^2$ bound for $\cT$
relies on an orthogonality argument based on the following result
(a rudimentary version of the orthogonality argument in the case of
two-sided
fold singularities is already in \cite{56}).

\proclaim{Lemma} For $\vectell=(l_0,\dots,l_m)$, with $0\le l_j\le l_0$
for $1\le j\le M$,
let $\cP_M(l)$ be the class of polynomials $\sum_{i=0}^M a_i h^i$ with
$2^{-l_i-2}\le |a_i|\le 2^{-\l_i+2}$ if $l_i<l_0$ and
$ |a_i|\le 2^{-l_0+2}$ if $l_i=l_0$.
Then there is a constant $C=C(M)$ and numbers $\nu_s,\mu_s$, $s=1,\dots,
10^M$
 so that

(i)
$$0\le \nu_1\le \mu_1\le \nu_2\le \mu_2\le \dots\le \nu_M\le \mu_M\le
1:=\nu_{M+1},$$

(ii) $$\nu_i\le \mu_i\le C\nu_i.$$

(iii) $$\Big|\sum_{i=1}^N a_i h^i\Big|\ge C^{-1} \max\{|a_j| |h|^j;
j=1,\dots M\}
\text{ if
$h\in [0,1]\setminus\bigcup_s[\nu_s,\mu_s].$}
$$
\endproclaim

Note that while $\mu_i$ and $\nu_i$ are close there may be `large' gaps between $\mu_i$ and 
$\nu_{i+1}$ for which the favorable lower bound (iii) holds.
The    elementary but somewhat lengthy   proof of the Lemma   based on
induction  is in \cite{66}.
A shorter and more elegant proof (of a
closely related inequality)  based on a
 compactness argument is due to Rychkov \cite{64}.

In order to descibe the orthogonality argument we  need some terminology.
Let $I$ be a subinterval  of $[0,1] $. We say that $\beta$ is a {\it
normalized cutoff function associated to $I$ } if
$\beta$ is supported in $I$ and $|\beta^{(j)}(t)|\le  |I|^{-j}$, for
$j=1,\dots,5$
and denote by $\fA(I)$ the set of all normalized cutoff functions associated
to $I$.

Fix $I$ and $\beta $ in $\fA(I)$; then we define another localization of
$\cT=T_{\la,\vectell}$  by
$$\cT[\beta] g(x)= \beta(x_1) \cT g(x).
\tag 3.15
$$

It follows quickly   from the definition and the property
$\nu_s\le \mu_s\le C\nu_s$   that
$$
\sup\Sb \widetilde I\\|\widetilde I|=\mu_s\endSb
\sup_{\widetilde \beta\in \fA(\Itil)}
%\sup\Sb \beta: \text{ norm. cutoff}\\
%\text{asso. to $I$}\endSb
\big \|\cT[\widetilde\beta]\big \|
\lc
\sup\Sb I\\|I|=\nu_s\endSb
\sup_{\beta\in \fA(I)}
%\sup\Sb\widetilde \beta: \text{norm. cutoff}\\
%\text{asso. to $\tilI$}\endSb
\big \|\cT[ \beta]\big \|.
\tag 3.16
$$
This is because for any interval $\widetilde I$ of length $\mu_s$ a function 
$\beta\in \fA(\widetilde I)$ can be written as a sum of a bounded number of functions associated to subintervals of length $\nu_s$.

We have to prove that also
$$
\sup\Sb \Itil\\|\Itil|=\nu_s\endSb
\sup_{\beta\in \fA(\Itil)}
%\sup\Sb\beta: \text{norm. cutoff}\\
%\text{asso. to $I$}\endSb
\big \|\cT[\widetilde\beta]\big \|
\lc
\sup\Sb I\\|I|=\mu_{s-1}\endSb
\sup_{\beta\in \fA(I)}
%\sup\Sb\widetilde \beta: \text{norm. cutoff}\\
%\text{asso. to $\tilI$}\endSb
\big \|\cT[\beta]\big \| + 2^{l_0(\frac 12+\eps)}\la^{-1}.
\tag 3.17
$$
and
$$\sup\Sb I\\|I|=\nu_1/8\endSb
\sup_{\beta\in \fA(I)}
\big \|\cT[\beta]\big \|\lc
 2^{l_0(\frac 12+\eps)}\la^{-1}.
\tag 3.18
$$

By the above remark inequality (3.17) is  obvious if $\mu_{s-1}\approx \nu_s$. Thus let us assume that
$\mu_{s-1}\le 2^{-100 M}\nu_s$ and fix $\widetilde \beta\in \fA(\tilI)$,
$|\tilI|=\mu_{s-1}$.

One uses the Cotlar-Stein Lemma
in the form
$$\Big\|\sum A_j\Big\|\lc \Big[\sum_{n=-\infty}^\infty \sup_j\|A^*_{j+n}
A_j\|^{\theta}\Big]^{1/2}
\Big[\sum_{n=-\infty}^\infty \sup_j\|A_{j+n} A_j^*\|^{1-\theta}\Big]^{1/2},
\tag 3.19
$$
for a (finite) sum of operators $\sum_{j} A_j$ on a Hilbert space.
(See \cite{73, ch. VII.2}; as pointed out
in \cite{7} and elsewhere, the version (3.19)
follows by a slight modification of the
standard  proof).

Now if
$J$ is an interval of length $\nu_s/8$ and
$\widetilde \beta\in \fA(J)$ then  we split
$\widetilde \beta=\sum_n\beta_n$ where for a fixed absolute constant $C$ the
function
$C^{-1}\beta_n$ belongs to $\fA(I_n)$  and the $I_n$  are intervals
 of length $\mu_{s-1}$;  $I_n$ and $I_{n'}$ are disjoint
if $|n-n'|>3$
 and the sum extends over no more than $O(\nu_s/\mu_{s-1})$ terms and thus over no more than
$O(2^{l_0})$ terms.

Now  let  $|I_n|=|I_{n'}|\approx \mu_{s-1}$ and $\dist (I_n,I_{n'})\approx
|n-n'||I|$ and assume 
$|n-n'||I| \le \nu_s/8$. Let $\beta_n$, $\beta_{n'}$ be  normalized cutoff
functions associated to $I_n$, $I_{n'}$. Then
$$\|\cT[\beta_n]^*\cT[\beta_{n'}]\|=0 \qquad \text{ if } |n-n'|>3
\tag 3.20
$$
by the disjointness of the intervals $I_n$, $I_n'$.
The crucial estimate is
$$\|\cT[\beta_n]\cT[\beta_{n'}]^*\|\lc
|n-n'|^{-1} 2^{l_0(1+\eps)}\la^{-2}
 \qquad \text{ if } |n-n'|>3, \ |n-n'|\le
\frac{\nu_s}{8\mu_{s-1}}. \tag 3.21
$$

(3.20/21) allows us to apply  (3.19) with $\theta=0$ (the standard version
does not apply,
as is erroneously quoted in \cite{66}).
This yields the bound
$$\|\cT[\widetilde \beta]\|\le C_\eps\big[ \sup_n\|\cT[\beta_n]\|+
2^{l(1+\eps)/2}\la^{-1} \sum_{n=3}^{2^{l_0}}n^{-1}\big]
$$
and thus (3.17).

To see (3.21) one examines the kernel $K$ of
 $\cT[\beta_n]\cT[\beta_{n'}]^*$
which is given by
$$
K(x,w)=\beta_{n}(x_1)\beta_{n'}(w_1)
\int e^{-\ic\la\tau(S(w,y_1)-S(x,y_1))}
b(x,w,y_1,\tau)
% \text{\it symbol }
 dy_1d\tau
\tag 3.22
$$
and by definition of $\mu_{s-1}$, $\nu_s$, $I_n$, $I_{n'}$ and the above
Lemma  we have
$$
|S_{y_1}(x,y_1)-S_{y_1}(w,y_1)|\gc 2^{-l_0}|x_1-y_1|- O(S(x,y_1)-S(w,y_1))
\tag 3.23
$$
To analyze the kernel $K$ and prove (3.21) by Schur's test
 one integrates by parts once in
$y_1$
and then many times in $\tau$, for the somewhat lengthy details
see \cite{66}, \cite{67}.
Analogous arguments also apply to the estimation of
$\cT[\beta]\cT[\beta]^*$ when
$\beta$ is associated to an interval of length $\ll\nu_1$, this gives
(3.18).

\subheading{Remarks}

{\bf 3.5.1.} Phong and Stein, in the remarkable paper  \cite{58},
 proved sharp $L^2$ decay estimates for oscillatory integral
operators with kernel $e^{i\la s(x,y)}\chi(x,y)$  in one dimensions,
where
$s$ is {\it real analytic}.
From their result and standard arguments one gets an improved result for the
generalized Radon
transform in the special semi-translation invariant case where  the curves
in $\Bbb R^2$ are
given by
$$y_2=x_2+s(x_1,y_1).\tag  3.24
$$
Namely, if $\cR f(x)=\int f(y_1, x_2+s(x_1,y_1)) \chi(x,y_1) dy_1$
then the endpoint $L^2\to L^2_\alpha$ estimate in Theorem 3.1 holds true.
An only  slightly weaker result for the case $s\in C^\infty$ has been obtained
by Rychkov \cite{64}. For related work  see also some   recent papers by
Greenblatt \cite{23}, \cite{24}.

{\bf 3.5.2.}
It is  not known exactly
which endpoint bounds hold in the general case of Theorem 3.3.
As an easy case  the $L^p\to L^p_{1/p}$ estimate holds if $p>n$ and
a type $(0, n-2)$ condition is satisfied (in the terminology of Theorem
3.3).
A similar statement for $1<p<n/(n-1)$ is obtained for type $(n-2,0)$ conditions
by passing to the adjoint operator.

The interpolation idea (3.8-9) is not limited
to conormal operators. Using variants of this method, sharp $L^p$ estimates
for Fourier integral operators in the nondegenerate case (\cite{68})
 were extended to  certain classes
with one- or two-sided fold singularities (\cite{70}, \cite{16}).
For other $L^p$ Sobolev  endpoint bounds in
 special cases
see \cite{74}, \cite{66},  \cite{57}, \cite{80}.

{\bf 3.5.3.}
Some endpoint inequalities in Theorem 3.3 fail:
M. Christ \cite{9} showed that  the convolution with a compactly supported
density on $(t,t^n)$ fails to map $L^n\to L^n_{1/n}$.
The best possible substitute is an $L^{n,2}\to L^n_{1/n}$ estimate
in \cite{69}; here $L^{n,2}$ is the Lorentz space.

{\bf 3.5.4.} Interpolation of the bounds in Theorem 3.3
with trivial $L^1\to L^\infty$ bounds (with loss of one derivative)
yields  almost sharp $L^p\to L^q$ bounds (\cite{66}, \cite{67}).
Endpoint estimates  for the case of two-sided finite type conditions
are in \cite{1}. For endpoint $L^p\to L^q$
estimates in the case (3.24),  with real-analytic $s$,
see
\cite{57}, \cite{79}, \cite{42}.

{\bf 3.5.5.} It would be desirable to obtain almost sharp $L^2$
versions such as Theorem 3.1 for more general oscillatory integral
operators with a corank one assumption.  Sharp endpoint $L^2$ results
 where one projection is a Whitney fold (type $1$)
and the other projection satisfies a finite type condition
are due to Comech \cite{15}.

{\bf 3.5.6.} Interesting bounds for the semi-translation
invariant case (3.24)
where only lower bounds on $s_{xy}$ (or higher derivatives) are assumed
were obtained by Carbery, Christ and Wright \cite{6}.
Related is the work by Phong, Stein and Sturm (\cite{60}, \cite{62},
\cite{61}), with
important contributions
concerning the stability of estimates.

\head{\bf 4. Operators with one-sided finite type conditions}\endhead

We now discuss operators of the  form (1.1) and assume that one of
the projections, $\pi_L$, is of type $\le r$ but make no assumption on the
other projection,
$\pi_R$. The role of the projections can be interchanged
by passing to the adjoint operator.

\proclaim{4.1.Theorem \cite{25},\cite{26},\cite{28}}
Suppose $\pi_L$ is of corank $\le 1$ and type $\le r$, and suppose
 that $\det d\pi_L$ vanishes simply.
If $r\in\{1,2,3\}$ then
$$
\|T_\la\|_{L^2\to L^2}\lc \la^{-(d-1)/2-1/(2r+2)}
\tag 4.1
$$
\endproclaim

It is conjectured that  this bound also holds for $r>3$.
The estimate  (4.1) is sharp in cases where the other projection
exhibits maximal degeneracy. In fact if $\pi_L$ is a fold and $\pi_R$ is a
blowdown then
more information is available such as a rather precise
 description of the kernel of $T_\la T_\la^*$, \cf.
Greenleaf  and Uhlmann \cite{32}, \cite{33}.
Applications include  the restricted X-ray transform in three dimensions
for the case where the line complexes are admissible in the sense of Gelfand
(\cite{21}, \cite{30}, \cite{34}); for an
early construction and application of a Fourier integral operator
 with this structure see also \cite{43}.

In the discussion that  follows we shall
replace the assumption that
$\det d\pi_L$ vanishes simply
(i.e., $\nabla_{x,z} \det\pi_L\neq 0$) by the more restrictive
assumption
$$
\nabla_z\det \pi_L\neq 0.
\tag 4.2
$$
In the case $r=1$ this is automatically satisfied, and it is shown in
\cite{26}, \cite{28} that in the cases $r=2$ and $r=3$ one can apply
canonical transformations to reduce matters to this situation.
For the oscillatory integral operators coming from the restricted X-ray
transform
for well-curved line complexes,
the condition
(4.2) is certainly satisfied.
We shall show that for general $r$
the estimate (4.1)  is a consequence of sharp estimates  for
oscillatory integral operators satisfying
 two-sided finite type conditions of order $r-1$, in $d-1$ dimensions.
The argument is closely related to Strichartz estimates
and can also be used to derive
$L^2\to L^q$  estimates
(an early version  can be found in Oberlin \cite{48}).

We shall now outline this argument.
After initial changes
of variables in $x$ and $z$ separately  we may assume that
$$
\aligned
&\Phi_{x'z'}(0,z)= I_{d-1}, \qquad
\Phi_{x'z_d}(0,z)= 0,
\\
\qquad
&\Phi_{x'z'}(x,0)= I_{d-1},\qquad
\Phi_{x_dz'}(x,0)= 0;
\endaligned
\tag 4.3.1
$$
moreover by our  assumption on the type  we may assume that
$$\Phi_{x_dz_d^{r+1}}(0,0)\neq 0
\tag 4.3.2
$$ and that
$\Phi_{x_dz_d^j}(x,z)$ is small for $j\le r$.
We may assume that the amplitudes are supported where $|x|+|z|\le \eps_0\ll
1$.

We form the operator $T_\la T_\la^*$ and write
$$T_\la T_\la^*f(x'x_d)=\int \cK^{x_d,y_d}[f(\cdot, y_d)](x') dy_d
\tag 4.4$$
where the kernel of
 $\cK^{x_dy_d}$  is given by
$$
K^{x_dy_d}(x',y')=\int
e^{\ic\la[\Phi(x',x_d,z)-\Phi(y',y_d,z)]}\sigma(x,z)\overline
{\sigma(y,z)} dz.
$$
We split
$K^{x_dy_d}=
H^{x_dy_d}+R^{x_dy_d}$ where
$H^{x_dy_d}(x',y')$ vanishes when
 $|x_d-y_d|\le \la^{-1}$ and when
$|x'-y'|\gc \eps|x_d-y_d|$, for some  $\eps$ with $\eps_0\ll \eps\ll 1$.

Notice that  by (4.3.1)
$$\Phi_{x'}(x',x_d,z)-\Phi_{x'}(y',y_d,z)=x'-y'+O(\eps_0|x_d-y_d|) $$
and by an integration by parts argument we
 get
$$|R^{x_dy_d}(x',y')|\le C_N (1+\la|x'-y'|)^{-N}
$$
for any $N$, in the relevant range $|x'-y'|\gc \eps|x_d-y_d|$.
Thus the corresponding operator $\cR^{x_dy_d}$ is bounded on $L^2(\Bbb
R^{d-1})$ and satisfies
$$\|\cR^{x_dy_d}\|_{L^2\to L^2}\le C_N' \la^{-d+1}(1+\la|x_d-y_d|)^{-N+d-1},
\tag 4.5$$
for any $N$.
For the main contribution $\cH^{x_dy_d}$ we are aiming for the estimate
$$
\|\cH^{x_dy_d}\|_{L^2\to L^2}
\lc \la^{1-d-\frac 1{r+1}} |x_d-y_d|^{-\frac 1{r+1}}.
\tag 4.6
$$
From (4.4), (4.5), (4.6)  and the $L^2(\Bbb R)$ boundedness of the
operator
with kernel $|x_d-y_d|^{-1/(r+1)}\chi_{[-1,1]}(x_d-y_d)$
the bound (4.1) follows  in a straightforward way.

Now observe that the operator $\cH^{x_dy_d}$ is local on cubes of diameter
$\approx |x_d-y_d|$ and we can use a trivial orthogonality argument to put
the localizations to cubes together. For a single cube we may then apply
a rescaling argument.
Specifically, let
$c\in \Bbb R^d$ and define
$$\widetilde H_c^{x_dy_d}(u,v)=H^{x_dy_d}(c+u|x_d-y_d|, c+v|x_d-y_d|).$$
Then for the corresponding operators we have
$$\|\cH^{x_dy_d}\|_{L^2\to L^2}\lc|x_d-y_d|^{d-1} \sup_c\|\widetilde {\Cal
H}_c^{x_dy_d}
\|_{L^2\to L^2}.
\tag 4.7
$$
Note that $\widetilde {\Cal H}_c^{x_dy_d}$ does not vanish only for small $c$.
A calculation shows that
the kernel of $\widetilde{ \Cal H}_c^{x_dy_d}$ is given by an oscillatory
integral
$$
\int e^{\ic\mu \Psi^\pm(u,v, z;c,\alpha,y_d)} b(u,v,z; c,x_d,y_d) dz; \qquad
\alpha=|x_d-y_d|,\,
\mu=\la|x_d-y_d|,
\tag 4.8
$$
with small parameters $c,\alpha=|x_d-y_d|,y_d$, and the phase function
is given by
$$
\Psi^\pm(u,v,z; \alpha,y_d,c)= \inn{u-v}{\Phi_{x'}(0,z)}\pm \Phi_{x_d}(0,z)
+\rho^{\pm}(u,v,z; \alpha,y_d,c).
$$
Here
the choice of $\Psi^{+}$ is taken if $x_d>y_d$ and
$\Psi^-$ is taken if $x_d<y_d$;
for the  error we have  $\rho^\pm= O(\alpha(|y_d|+c))$ in the $C^\infty$
topology.
Observe in particular that  for $\alpha=0$ we get essentially the
localization of
a translation invariant operator.

We now examine the canonical relation associated to the
 oscillatory integral, when $\alpha=0$.
In view of (4.3.1) the critical set $\{\nabla_z\Psi^\pm =0\}$
 for the phase function  at $\alpha=0$ is given by
$\{(u,v,z): v=u+g(z), \Phi_{x_dz_d}(0,z)=0\}$
for suitable $g(z)$;  in view of (4.2) this defines a smooth manifold.
Consequently the canonical relation
$$\cC_{\Psi^\pm}\Big|_{\alpha=0}=\{(u, \Psi^\pm_u,  v,
\Psi^\pm_v):\Psi^\pm_z=0\}$$ is a smooth manifold. By (4.2) we may assume
(after performing a rotation)  that
$\Phi_{x_dz_dz_1}\neq 0$ and then solve the equation $\Phi_{x_dz_d}(0,z)=0$
near the origin  in terms of a function $z_1=\widetilde z_1(z'', z_d)$.
The projection $\pi_L$ is  given by
$$(u,z'',z_d)\to (u,\Phi_{x'}(0,z_1^\pm(z'',z_d), z'',z_d))$$
and $\partial /\partial z_d$ is a kernel field for
$\pi_L$.
Implicit differentiation reveals that
$\partial_{z_d}^kz_1^\pm - \Phi_{z_dx_dz_1}^{-1}\Phi_{x_dz_d^{k+1}}$
 belongs to the ideal generated by $\Phi_{x_dz_d^{j}}$, $j\le k$ and
thus, by our assumption (4.3.2) we see that $\pi_L$ is  of type
$\le r-1$. The same holds true for $\pi_R$, by symmetry considerations.
Although we have verified these conditions  for $\alpha=0$ they remain true
for small $\alpha$ since Morin singularities are stable under small
perturbations.

We now discuss estimates for the oscillatory integral operator
$S^\pm _\mu$ whose kernel is given by (4.8) (we suppress
the dependence on $ c,\alpha, y_d$.)
The number of frequency variables  is $N=d$ and thus
we can expect the uniform bound
$$
\|S^\pm _\mu\|_{L^2\to L^2}\lc \mu^{-\frac{d-2}2-\frac 1{r+1} -\frac{d}2}
\tag 4.9
$$
for small $\alpha$.
Indeed,
the case $\alpha=0$ of (4.9) is easy to verify;
 because of the
translation invariance we may apply Fourier transform arguments together
with the method of stationary phase and van der Corput's lemma.
Given (4.9)
we obtain from (4.7) and from (4.9) with $\mu=\la|x_d-y_d|$ that
$$
\|\cH^{x_dy_d}\|_{L^2\to L^2}\lc |x_d-y_d|^{d-1}
\mu^{-(d-1)-\frac 1{r+1} }\lc \la^{-(d-1)-\frac 1{r+1}}
|x_d-y_d|^{-\frac 1{r+1}}.
$$

Of course the Fourier transform argument   does not extend to the
case where $\alpha$ is  merely small.
However  if $r=1$ the estimate follows from (1.6) (with
$d$ replaced by $d-1$ and
 $N=d$)  since then $\cC_{\Psi^\pm}$ is
a local canonical graph.
Similarly, if $r=2$ then the canonical relation $\cC_{\Psi^\pm}$ projects
with two-sided fold
singularities so that the desired estimate follows from known estimates for
this situation
(see the pioneering paper by Melrose and Taylor \cite{44},
and also \cite{53}, \cite{19}, \cite{27}).
For the case $r=3$, inequality (4.6) follows from a recent result by
the authors \cite{28} discussed in the next section,
plus the reduction outlined in \S1.2. The case $r\ge 4$ is
currently open.

\remark{\bf  Remarks}

{\bf 4.2.1.} The argument  above can also be used to prove $L^2\to L^q$
estimates (see \cite{48}, \cite{25}, \cite{26}).
Assume $r=1$ and thus  assume that $\pi_L: \cC\to T^*\Om_L$ projects  with
Whitney folds. Then a stationary phase argument gives that
$$\|\cK^{x_dy_d}\|_{L^1\to L^\infty} \lc (1+\la|x_d-y_d|)^{-1/2}
\tag 4.10$$
and interpolation with (4.5-6) yields $L^{q'}\to L^{q}$ estimates for
$\cK^{x_dy_d}$ and then $L^2\to L^q$ bounds for $T_\la$.
The result \cite{25} is
$$
\|T_\la\|_{L^2\to L^q}\lc \la^{-d/q}, \qquad 4\le q\le \infty.
\tag 4.11
$$

The estimate (4.10) may be improved under the presence of some curvature
assumption.
Assume that the projection of the fold surface $S_1(\pi_L)$ to $\Omega_L$
is a submersion, then for each $x\in \Omega_L$
the projection of $S_1(\pi_L)$ to the fibers is a hypersurface $\Sigma_x$ in
$T_x^*\Omega_L$. Suppose that for every $x$
this hypersurface has $l$ nonvanishing principal curvatures
(this assumption is reminiscent of the
so-called cinematic curvature hypothesis in \cite{46}). Then
(4.10) can be replaced by
$$\|\cK^{x_dy_d}\|_{L^1\to L^\infty} \lc (1+\la|x_d-y_d|)^{-(l+1)/2}
$$
and
(4.11)
holds true for a larger range of exponents, namely
$$\|T_\la\|_{L^2\to L^q}\lc \la^{-d/q}, \qquad \frac{2l+4}{l+1}
\le q\le \infty.
\tag 4.12$$

The version of this estimate for Fourier integral operators
\cite{25}, with $l=1$, yields Oberlin's sharp
$L^p\to L^q$  estimates \cite{48} for the averaging operator (2.7)
 in
three dimension (assuming that $\Gamma$ is nondegenerate), as well as
 variable coefficient perturbations. It also yields sharp results for
certain
convolution operators associated to curves
 on the Heisenberg group (\cite{65}, see  \S7.3 below)
and for   estimates for restricted X-ray transforms
associated to well curved line complexes in  $\Bbb R^3$ (\cite{25}).

In  dimensions $d>3$ the method yields $L^2\to L^q$ bounds (\cite{26})
which should be considered as partial results, since in most
interesting cases
the endpoint $L^p\to L^q$ estimates do not involve the exponent $2$.

{\bf 4.2.2.}
The analogy with the cinematic curvature hypothesis has been exploited by
Oberlin, Smith and Sogge \cite{52}
to prove nontrivial $L^4\to L^4_\alpha$
estimates for translation invariant operators associated to
nondegenerate curves in $\Bbb R^3$. Here it is crucial to apply
a  square function estimate due to Bourgain \cite{3} that he used
 in proving bounds  for cone multipliers.
The article \cite{51} contains an interesting counterexample
for the failure of
$L^p\to L^{p}_{1/p-\eps}$ estimates when $p<4$.

{\bf 4.2.3.} Techniques of
 oscillatory integrals have been used by Oberlin
\cite{49} to obtain essentially sharp $L^p\to L^q$ estimates for the
operator (2.7) in four dimension, see also  \cite{29} for a
related argument for the
restricted X-ray transform in four dimensions, in the rigid case (2.8).

{\bf 4.2.4.} More recently, a  powerful combinatorial
method was developed by Christ \cite{10} who proved essentially sharp
$L^p\to L^q$
estimates for the translation invariant model operator (2.7)
in all dimensions
(for nondegenerate $\Gamma$).
$L^p\to L^q$ bounds for the X-ray transform  in higher dimensions,
in the model case (2.8), have been obtained by Burak-Erdo\u gan and Christ
(\cite{4},
 \cite{5}); these papers  contain even stronger
 mixed norm estimates.
Christ's combinatorial
 method has been  further developed by Tao and Wright \cite{75} who
obtained
almost sharp $L^p\to L^q$  estimates for variable coefficient analogues.

\endremark

\head{\bf 5. Two-sided type two  singularities}\endhead

We consider again the operator (1.1) and discuss the
 proof of the following
result mentioned in the last section.

\proclaim{5.1. Theorem \cite{28}} Suppose that both
$\pi_L$ and $\pi_R$ are of type $\le 2$. Then for $\la\ge 1$
$$\|T_\la\|_{L^2\to L^2}=O( \la^{-(d-1)/2-1/4}).$$
\endproclaim

A slightly weaker version of this result is due to Comech and Cuccagna
\cite{17}
who obtained the bound
$\|T^\la\|\le C_\eps \la^{-(d-1)/2-1/4+\eps}$ for $\eps>0$.

The proof of the endpoint estimate 
is  based on various localizations and almost orthogonality
arguments.
As in \S2 we start with localizing the determinant of
$d\pi_{L/R}$ and its
derivatives with respect to a kernel vector field.
The  form (5.2) below of this first  decomposition can already be found in
\cite{15}, \cite{17}.

We assume that the amplitude is supported near the origin and assume that
(4.3.1) holds. Let $\Phi^{z'x'}=\Phi_{x'z'}^{-1}$,
$\Phi^{x'z'}=\Phi_{z'x'}^{-1}$; then
 kernel vector fields for the projections $\pi_L$ are given by
$$\aligned
V_R&=\partial_{x_d}-\Phi_{x_dz'}\Phi^{z'x'}\partial_{x'},
\\
V_L&=\partial_{z_d}-\Phi_{z_dx'}\Phi^{x'z'}\partial_{z'},
\endaligned
\tag 5.1
$$
respectively. Also let $h(x,z)=\det \Phi_{xz}$ and by  the type two
assumption we can assume that
$|V_L^2 h|$,
$|V_R^2 h|$ are bounded below.
Emphasizing the amplitude in (1.1) we write $T_\la[\sigma]$ for the operator
$T_\la$ and will introduce various decompositions of the amplitude.

Let $\beta_0\in C^\infty(\Bbb R)$ be an even function supported in $(-1,1)$,
and equal to one in $(-1/2,1/2)$ and for
$j\ge 1$ let
$\beta_j(s)=
\beta_0(2^{-j}s)-\beta_0(2^{-j+1}s)$.
Denote by $\ell_0$
 that is  the largest integer $\ell$
 so that $2^{\ell}\le \la^{1/2}$ (we assume that $\la$ is large).
Define
$$
\aligned
\si_{j,k,l}(x,z)&= \si(x,z) \beta_1(2^l h(x,z))\beta_j( 2^{l/2}V_Rh(x,z))
\beta_k(2^{l/2}V_Lh(x,z))
\\
\si^0_{j,k,\ell_0}(x,z)&=
 \si(x,z) \beta_0(2^{\ell_0} h(x,z))\beta_j( 2^{\ell_0/2}V_Rh(x,z))
\beta_k(2^{\ell_0/2}V_Lh(x,z));
\endaligned
\tag 5.2
$$
thus if $j,k>0$ then
$|h|\approx 2^{-l}$,
$|V_L h|\approx 2^{k-l/2}$,
$|V_R h|\approx 2^{j-l/2}$
 on the support of $\si_{j,k,l}$.

It is not hard to see that the estimate of Theorem 5.1 follows from

\proclaim{ 5.2. Proposition}
We have the following bounds:

(i) For $0<l<\ell_0=[\log_2(\sqrt\la)]$
$$
\|T_\la[\si_{j,k,l}]\|_{L^2\to L^2}\lc
\la^{-(d-1)/2}\min\big \{ 2^{l/2}
\la^{-1/2}; 2^{-(l+j+k)/2}
\big\}.
\tag 5.3
$$

(ii)
$$
\|T_\la[\si^0_{j,k,\ell_0}]\|_{L^2\to L^2}\lc
\la^{-(d-1)/2-1/4} 2^{-(j+k)/2}.
\tag 5.4
$$
\endproclaim

We shall only discuss (5.3) as (5.4) is proved similarly. In what follows
$j,k,l$ will be fixed and we shall discuss the main case
where $0<k\le j\le l/2$,
$2^l\le \la^{1/2}$.
As in the argument in \S2 standard $T^*T$ arguments do not work
and further localizations and almost orthogonality arguments are needed.
These are less straightforward in the higher dimensional situation
considered here, and the amplitudes will be localized to nonisotropic boxes
of various sides depending on the geometry of the kernel vector fields.

For $P=(x^0,z^0)\in\Omega_L\times \Omega_R $  let
$ a_P=\big(
- \Phi^{x'z'}(P) \Phi_{z'x_d}(P),1\big)$ and
$b_P=(-\Phi^{z'x'}(P)\Phi_{x'z_d}(P),1)$ so that
 $V_L=\inn{a_P}{\partial_x}$,
 $V_R=\inn{b_P}{\partial_z}$.
Let  $\pi_{a_P}^\perp$, $\pi_{b_P}^\perp$ be the orthogonal projections  to
the orthogonal complement of $\bbR a_P$ in $T_{x^0}\Omega_L$ and
$\bbR b_P$ in $T_{z^0}\bbR^d$,
respectively. Suppose  $0<\gamma_1\le \gamma_2\ll 1$ and
$0<\delta_1\le\delta_2\ll 1 $   and
let
$$
B_P(\gamma_1,\gamma_2,\delta_1,\delta_2)$$
denote the box of all $(x,z)$ for which
$|\pi_{a_P}^\perp(x-x^0)|\le \gamma_1$, $|\inn{x-x^0}{a_P}|\le \gamma_2$,
$|\pi_{b_P}^\perp(z-z^0)|\le \delta_1$, $|\inn{z-z^0}{b_P}|\le \delta_2$.
We always assume
$$\gamma_1\le \gamma_2, \ \delta_1\le \delta_2 \tag 5.5$$
 We say that $\chi\in C^\infty_0$ is a {\it  normalized cutoff function
associated to $B_P(\gamma_1,\gamma_2,\delta_1,\delta_2)$} if it is
supported in
$B_P(\gamma_1,\gamma_2,\delta_1,\delta_2)$ and satisfies the (natural)
estimates
$$
|(\pi_a^\perp\nabla_x)^{m_L}
\inn{a}{\nabla_x}^{n_L}
(\pi_b^\perp\nabla_z)^{m_R}
\inn{b}{\nabla_z}^{n_R}
\chi(x,z)|\le  \gamma_1^{-m_L}\gamma_2^{-n_L}
 \delta_1^{-m_R}\delta_2^{-n_R}
$$
whenever  $m_L+n_L\le 10d$,  $m_R+n_L\le 10d$.

We denote by $\fA_P(\gamma_1,\gamma_2,\delta_1,\delta_2)$ the class of
all normalized cutoff functions associated to
$B_P(\gamma_1,\gamma_2,\delta_1,\delta_2)$.

Suppose that
$(\gamma_1,\gamma_2,\delta_1,\delta_2)=
(\eps 2^{-l}, \eps 2^{-j-l/2},\eps 2^{-l}, \eps 2^{-k-l/2})$. It turns
out  that
% in view  of the geometry of the maps $\pi_L$ and $\pi_R$
$h=\det\Phi_{xz}$ changes only by $O(\eps 2^{-l})$ in the box
$B_P(\gamma_1,\gamma_2,\delta_1,\delta_2)$ but is in size comparable to
$2^{-l}$. This enables one to apply a $TT^*$ argument and one obtains
the correct bound $ O(2^{l/2}\la^{-d/2})$
for the operator norm of $T_\la[\chi\sigma]$ assuming that $\chi\in
\fA_P(\gamma_1,\gamma_2,\delta_1,\delta_2)$ for some fixed $P$.
This step had already been carried out by Comech and Cuccagna
\cite{17}.
Let
$$\cA_P(\gamma_1,\gamma_2,\delta_1,\delta_2):=\sup\big\{
\big\|T_\la[\chi \si_{j,k,l}]\big\|: \ \chi\in
\fA_P(\gamma_1,\gamma_2,\delta_1,\delta_2)\big\}
$$
then,
for $2^l\le \la^{1/2}$,
$$
\sup_P \, \cA_P(2^{-l}, 2^{-j-l/2}, 2^{-l}, 2^{-k-l/2})
\lc 2^{l/2}\la^{-d/2}.
\tag 5.6
$$
If one  uses that $|V_R h|\approx 2^{j-l/2}$,
$|V_L h|\approx 2^{k-l/2}$ one also
gets
$$
\sup_P \, \cA_P(2^{-l}, 2^{-j-l/2}, 2^{-l}, 2^{-k-l/2})
\lc 2^{-(l+j+k)/2}\la^{-(d-1)/2}.
\tag 5.7
$$
Initially one obtains these estimates for boxes of size
$(\eps 2^{-l}, \eps 2^{-j-l/2},\eps 2^{-l}, \eps 2^{-k-l/2})$
but the $\eps $ may be removed since we can decompose any
$B_P(\gamma_1,\gamma_2,\delta_1,\delta_2)$
into no more than
$O(\eps^{-2d})$ boxes
of dimensions
 $(\eps \gamma_1,\eps \gamma_2,\eps \delta_1,\eps \delta_2)$.
From this one deduces
$$\cA_P(\gamma_1,\gamma_2,\delta_1,\delta_2)\le C_\eps
\sup_Q\cA_Q(\eps \gamma_1,\eps\gamma_2,\eps\delta_1,\eps\delta_2).
\tag 5.8
$$

In order to put the localized pieces together we need some orthogonality
arguments. For the sharp result we need to prove various
 inequalities of the
form
 $$\multline
\sup_P\,\cA_P(\gamma_{1,\lrg}, \gamma_{2,\lrg},
\delta_{1,\lrg}, \delta_{2,\lrg})
\\
\lc
\sup_Q\cA_Q(\gamma_{1,\sm}, \gamma_{2,\sm},
\delta_{1,\sm}, \delta_{2,\sm})
+E(j,k,l)
\endmultline
\tag 5.9
$$
where the error term satisfies
$$E(j,k,l)\lc \la^{-(d-1)/2} \min\{ 2^{l/2}\la^{-1}, 2^{-(l+j+k)/2}\}
\tag 5.10
$$
or a better estimate.

In the argument it is crucial that we assume
$$\min \{\frac{\gamma_{1,\sm}}{\gamma_{2,\sm}},
\frac{\delta_{1,\sm}}{\delta_{2,\sm}}\}\gc
\max\{\gamma_{2,\lrg},\delta_{2,\lrg}\}
\tag 5.11
$$ since from (5.11) one can see that the orientation of small boxes
$B_Q(\gamma_\sm,\delta_\sm)$ does not  significantly change if
$Q$ varies in the large box $ B_P(\gamma_{\lrg},\delta_{\lrg}) $.

\proclaim{5.3 Proposition}
Let $k\le j\le l/2, 2^l\le \la^{1/2}$.
There is $\eps>0$ (chosen independently of $k,j,l,\la$)
so that the   inequality (5.9) holds with the choices of

\roster
\item"{(i)}"
$$\aligned
(\gamma_{\lrg}, \delta_{\lrg})&=
(\eps 2^{j+k-l}, \eps 2^{k-l/2}, \eps 2^{j+k-l}, \eps 2^{k-l/2})
\\
(\gamma_{\sm}, \delta_{\sm})&=
(2^{-l}, 2^{-j-l/2}, 2^{-l}, 2^{-k-l/2}) ,
\endaligned
\tag 5.12
$$

\item"{(ii)}"
$$\aligned
(\gamma_{\lrg}, \delta_{\lrg})&=
(\eps 2^{j-l/2}, \eps 2^{j-l/2}, \eps 2^{k-l/2}, \eps 2^{k-l/2})
\\
(\gamma_{\sm}, \delta_{\sm})&=
(2^{j+k-l}, 2^{k-l/2}, 2^{j+k-l}, 2^{k-l/2}),
\endaligned
\tag 5.13$$

\item"{(iii)}"
$$\aligned
(\gamma_{\lrg}, \delta_{\lrg})&=
(\eps, \eps, \eps, \eps )
\\
(\gamma_{\sm}, \delta_{\sm})&=
(2^{j-l/2}, 2^{j-l/2}, 2^{k-l/2}, 2^{k-l/2}).
\endaligned
\tag 5.14
$$
\endroster
\endproclaim

A combination of these estimates (with 5.8)   yields the
desired bound (5.3); here the outline of the
argument is similar to the one given in \S3.
For each instance we are given a cutoff function $\zeta\in \fA_P
(\gamma_{\lrg}, \delta_{\lrg})$
and we decompose $$\zeta=\sum_{(X,Z)\in
\Bbb Z^d\times \Bbb Z^d} \zeta_{XZ}$$ where the
$\zeta_{XZ}$ is, up to a constant, a normalized cutoff function
associated to a box of dimensions $(\gamma_{\sm},\delta_{\sm})$; the various
boxes have bounded overlap, and comparable orientation. More precisely
if $(P,Q)$ is a reference point in the big box
$B_P(\gamma_{\lrg}, \delta_{\lrg})$ then  each of the small boxes is
 comparable to a box defined by the conditions
$|\pi_{a_P}^\perp(x-x_{X})|\le \gamma_1$, $|\inn{x-x_X}{a_P}|\le
\gamma_2$,
$|\pi_{b_P}^\perp(z-z_Z)|\le \delta_1$, $|\inn{z-z_Z}{b_P}|\le \delta_2$.

If $T_{XZ}$ denotes the operator $T_\la[\zeta_{XZ} \sigma_{j,k,l}]$ then in
each case we have to show that for large $N$
$$
\align
\|T_\la[\zeta_{XZ}] &(T_\la[\zeta_{X'Z'}])^*\|_{L^2\to L^2}
+
\|(T_\la[\zeta_{XZ}])^* T_\la[\zeta_{X'Z'}]\|_{L^2\to L^2}
\\
&\lc\la^{1-d} \min\{ 2^{l}\la^{-2}, 2^{-l-j-k}\} (|X-X'|+|Z-Z'|)^{-N}
\endalign
$$
if $|X-X'|+|Z-Z'|\gg 1$.

For the estimation in the case (5.12) it is crucial that in any fixed large
box
% $B_P(\eps 2^{j+k-l}, \eps 2^{k-l/2}, \eps 2^{j+k-l}, \eps 2^{k-l/2})$
$V_L h$ does
not change by more than $O(\eps 2^{k-l/2})$ and
thus is comparable to $2^{k-l/2}$ in the entire box; similarly
$V_R h$ is comparable to $2^{j-l/2}$ in the entire box. For the
orthogonality we use that
$\Phi_{x'z'}$ is close to the identity.
In the other extreme case (5.14) $V_R h$ and $V_L h$ change significantly
in the direction of kernel fields  and this can be exploited in the
orthogonality argument. (5.13) is an intermediate case.  
This description is
an oversimplification and we refer the reader to  \cite{28}
for the detailed discussion of each case.

\head
{\bf 6. Geometrical conditions on families of curves}
\endhead

We illustrate some of the results mentioned before by
relating conditions involving strong Morin singularities to
various conditions on vector fields and their commutators.

\subheading{6.1.
 Left and right  commutator conditions and strong Morin singularities}
%and relate conditions for strong cusps and finite type  for the
projections
%$\pi_L$ and $\pi_R$ to commutator conditions for
%$(1,0)$ and $(0,1)$  vector fields.
%
We first  look at an incidence relation $\cM$ with canonical relation
$\cC=N^*\cM$ as in
(1.12) and assume $\ell=d-1$ so that $\dim \cM=d+1$. As in \S3
\cite{67}, we have two distinguished
classes  of vector fields on $\cM$, namely vector fields of type
 $(1,0)$ which are also tangent to $\cM\cap (\Omega_L\times 0)$ and  vector
fields  of type
$(0,1)$ which are
tangent to $\cM\cap (\{0\}\times \Omega_R)$. Note that for each point $P$
the corresponding distinguished tangent spaces $T^{1,0}_P\cM$
and
 $T^{0,1}_P\cM$ are one-dimensional.
If $\Phi$ is the $\Bbb R^{d-1}$-valued defining function for
$\cM=\{\Phi(x,y)=0\}$ then
a nonvanishing $(1,0)$  vector field $X$ and 
a nonvanishing $(0,1)$  vector field $Y$ are given by
$$
X=
\sum_{j=1}^d
a_j(x,y)
\frac{\partial}{\partial x_j},
\quad
Y= \sum_{k=1}^d b_k(x,y)
 \frac{\partial}{\partial y_k}
\tag 6.1$$
where
$(-1)^{j-1} a_j(x,y)$ is the determinant
of the $(d-1)\times(d-1)$ matrix obtained from
the $(d-1)\times d$ matrix $\Phi_x'$ by omitting the
$j^{\text th} $
 column, and
$(-1)^{k-1}b_j(x,y)$ is the determinant of the  $(d-1)\times(d-1)$ matrix
obtained from
$\Phi_y'$ by omitting the $j^{\text th} $ column.

The canonical relation $N^*\cM'$ in (1.12)  can be identified with a
subbundle
 $T^{*,\perp}\cM$ of $T^*\cM$
whose fiber at $P\in\cM$ is
 the $\ell$-dimensional space of all linear
functionals in
$T^*_P\cM$ which annihilate vectors in $T^{1,0}_P\cM$ and vectors in
$T^{0,1}_P\cM$,
$$T^{*,\perp}_P\!\cM   = (T^{1,0}_P\cM\oplus T^{0,1}_P\cM)^\perp.$$
Concretely,
%The identification is via restricting linear forms in $T_P^*(\fX\times\fY)$
%to  tangent vectors in $T_P\cM$.
if
$\imath:\cM\to \Om_L\times \Om_R$ denotes the inclusion map and
$\imath^*$ the pullback of $\imath$  (or restriction operator)
acting on forms in $T^*(\fX\times\fY)$,  then
$$
T^{*,\perp}\cM=\{(P,
 \imath^*_P\la): \, (P,\la)\in\Cal C\}
$$

Finite type conditions can be formulated in terms of iterated commutators
of $(1,0)$ and $(0,1)$ vector fields (\cite{67}).
%
%For an explicit formula relating the rotational curvature (1.14) to
%$\sum\tau_i\inn{d_x\Phi^i-d_y\Phi^i}{[X,Y]}$ \cf. Lemma 2.3 in \cite{66}.
Here they are used to characterize the situation of strong
 Morin singularities (\cf. \S 2.4).
Let $x^0\in \Omega_L$, let $\cM_{x^0}=\{y\in \Om_R: (x^0,y) \in \cM\}$
 and let
$$
\fN_{L, x^0}:=\pi_L^{-1}(\{x^0\}\times T^*_{x_0}\Omega_L)=
\{(y,\la): y\in\cM_{x^0}, \, \la\in T^{*,\perp}_{(x^0,y)}\cM\}.
$$
Let $\pi_{L,x_0}$ the restriction of $\pi_L$ to $\fN_{L,x_0}$ as a
map to $T_{x_0}^*\Omega_L$,
then $\pi_L$ has strong Morin singularities if for fixed $x_0$ the map
$\pi_{L,x_0}$ has Morin singularities.

Similarly, if
 $y^0\in \Omega_L$, let $\cM^{y^0}=\{x\in \Om_L: (x,y^0) \in \cM\}$
then the adjoint operator $\cR^*$ is an integral operator along the
curves  $\cM^{y_0}$; now we  define
$\fN_{R, y^0}$ as the set of all
$(x,\la)$ where $x\in\cM_{y^0}$, $ \la\in T^{*,\perp}_{(x,y^0)}\cM$, and
$\pi_{R,y^0}:\fN_{R,y^0}\to T_{y^0}^*\Omega_R$
 is the restriction of the map $\pi_R$.

\proclaim{Proposition}

(a)  Let $x^0\in\Om_L$ and $y^0\in \cM_{x^0}$ and let $P=(x^0, y^0)$.
 The following statements are equivalent.

\roster
\item"{(i)}"  Near $P$, the only singularities of
$\pi_{L,x^0} $  are  $S_{1_{k},0}$ singularities, for $k\le d-2$.

\item"{(ii)}"
The vectors
$(\ad Y)^mX$, $m=1,\dots, d-1$ are linearly independent at $P$.
\endroster

(b)
 Let $y^0\in\Om_R$ and $x^0\in \cM^{y^0}$ and let $P=(x^0, y^0)$.
The following statements are equivalent.

\roster
\item"{(i)}"  Near $P$, the only singularities of
$\pi_{R,y^0} $  are  $S_{1_{k},0}$ singularities, for $k\le d-2$.

\item"{(ii)}"
The vectors
$(\ad X)^mY$, $m=1,\dots, d-1$ are linearly independent at $P$.
\endroster
\endproclaim

It suffices to verify statement (a).
There are coordinate systems
$x=(x',x_d)$ near $x^0$, vanishing at $x^0$
 and $y=(y',y_d)$ near $y^0$,  vanishing at  $y^0$ so that
near $P$ the manifold $\cM$ is given by $y'=S(x,y_d)$ with
$$
S(x,y_d)= x'+ x_d g(y_d)
+O(|x|^2)
%+x_d^2 b(x,y_d)+x_d \sum_{i=1}^{d-1}x_iA_i(x,y_d) +\sum_{i,j}x_ix_j
%R_{ij}(x,y_d)
$$
where $g(0)=0$.
%% $A_i$, $R_{ij}$  are smooth vector valued functions, with $g(0)=0$.

In these coordinates we compute the vector fields $X$ and $Y$ in (6.1) and
find
$$\aligned
(-1)^{d-1}
a_j&=  g_j(y_d)+O(|x|) , \qquad j=1,\dots, d-1,
\\
(-1)^{d-1}
a_d&= 1+O(|x|),
\endaligned
$$
and
$$\aligned
b_j
%= \frac{\partial S^j}{\partial y_d}
&= x_d\frac
{\partial g_j}{\partial y_d}+O(|x'|^2+|x'||x_d|), \qquad j=1,\dots, d-1,
\\
b_d&= 1.
\endaligned
$$
By induction one verifies that for $m=1,2,\dots$
$$(-1)^{d-1}(\ad Y)^mX= \sum_{j=1}^d v_j^m \frac{\partial}{\partial x_j}
+\sum_{j=1}^d w_j^m
\frac{\partial}{\partial y_j}
$$
where
$$\aligned
v_j^m&=  \frac{\partial^m g_j}{\partial y_d^m}+O(|x|) , \qquad j=1,\dots,
d-1,
\\
v_d^m&=O(|x|)
\endaligned
$$
and
$$\aligned
w_j^m&= - \frac{\partial^m g_j}{\partial y_d^m}+O(|x|) , \qquad j=1,\dots,
d-1,
\\
w_d^m&=O(|x|)
\endaligned
$$
Consequently we see that the linear independence of the vector fields
$(\ad Y)^m X$ at $P$ is equivalent with the linear independence of
$ \partial^m g_j/(\partial y_d^m)$ at $y_d=0$.

Next, the map $\pi_{L,x^0}:\fN_{L,x^0}\to T^*_{x^0}\Omega_L$
is in the above  coordinates given by
$$
(y_d, \tau)\mapsto \tau\cdot S_x(0,y_d)= (\tau_1,\dots, \tau_{d-1},
\sum_{i=1}^{d-1}\tau_i
g_i(y_d))
$$
and from (2.3-5) we see that the statement (i)
is also  equivalent with the linear independence of the vectors
$ \partial^m g_j/(\partial y_d^m)$ at $y_d=0$.

This proves the proposition.

\subheading{6.2.
Families of curves defined by exponentials of vector fields}
Let now $\{\gamma_t(\cdot)\}_{t\in I}$ be a one-parameter family of
diffeomorphisms of $\br^n$ which we  can also
consider as a family of parametrized curves,
$$t\mapsto \gamma_t(x):= \gamma(x,t).$$
We shall  assume that $x$ varies in an open set $\Omega$, the open parameter
interval $I$ is a small neighborhood of $0$ and
that
$\g_0 =Id$ and $\dot\g\ne 0$, where $\dot\g$
denotes $\frac{d}{dt}(\g_t)$. Thus
for each $x$, $t\mapsto \gamma(x,t)$
 defines a regular curve passing through $x$.
As in the article by Christ, Nagel, Stein and Wainger
\cite{11}, we may write such a family as
$$\gamma_t(x):=\gamma(x,t)=\exp(\sum_{i=1}^N t^i X_i)(x)\mod
O(t^{N+1})\tag 6.2$$
for some vector fields $X_1, X_2,...$, and $N\in \Bbb N$.
The generalized  Radon transform is now defined by
$$\Cal Rf(x)=\int f(\g(x,t))\chi(t) dt
$$
and  incidence relation $\cM$ is given by
$$\cM=\{(x,\gamma(x,t)):x\in\br^n, t\in\br\}\subset\br^n\times\br^n.
\tag 6.3$$

Besides using the projections $\pi_L$ and $\pi_R$, there are other
ways of describing what it means for the family
$\{\gamma_t(\cdot)\}$ to be
maximally nondegenerate, in either a one- or two-sided fashion.
One is given in terms
the structure of the pullback map with respect to the diffeomorphisms
$\g_t(\cdot)$, and another is given by
the linear independence of certain linear combinations of the vector
fields
$X_j$
and their iterated commutators. We formulate the conditions
on the right, with the analogous conditions on the left being easily
obtained
by symmetry.

\subheading{\bf 6.3. Strong Morin singularities  and pull-back conditions}
We  are working with (6.2) and  formulate the {\it  pullback condition}
$(P)_R$.
Form the curve
$$\Gamma_R(x,t)=\frac{d}{ds}\bigl(\gamma_{s+t}\circ\gamma_t^{-1}(x)\bigr)
\bigr|_{s=0},
\tag 6.4$$
so that $\Gamma_R(x,\cdot):\br\rta T_x\br^n$.
Let  $\Gamma_R^{(\nu)}(x,t)=(\partial/\partial t)^\nu \Gamma_R(x,t)$ for
$\nu=0,1,\dots$

\definition{Definition}
The family of curves $\{\gamma(x,\cdot)\}_{x\in \Omega}$
satisfies condition $(P)_R$ at $x$
if the vectors $\Gamma_R^{(\nu)}(x,0)$, $\nu=0,\dots, n-1$ are linearly
independent.
\enddefinition

Let $\cM$ be the incidence relation for our averaging operator.

\proclaim{Proposition} Let $c_0=(x_0,\xi_0,x_0,\eta_0)\in N^*\cM'$.
Then condition $(P)_R$ is satisfied at $x_0$ if and only if
 $\pi_R$ has only $S_{1_k,0}^+$ singularities at $c$, with $k\le d-2$.
\endproclaim

To see this, note that
 $\cM\subset \R^n\times\R^n$  is the image of the immersion
$(x,t)\mapsto
(x,\gamma(x,t))$. Thus  $(x,\xi;y,\eta)$ belongs to $N^*\cM'$ if and only if
$y=\gamma(x,t)$ for some $t\in\R$
and
$(D\Phi_{(x,t)})^*(\xi,-\eta)=(0,0)\in T^*_{(x,t)}\R^{n+1}$. This yields
$$N^*\cM'=\bigl\{(x,(D_x\g)^*(\eta);\g(x,t),\eta):x\in\R^n, t\in\R,
\eta\cdot
\dot\g(x,t)=0\bigr\}.$$
For each fixed $t$, let $y=\g_t(x)$, so that $x=\g^{-1}_t(y)$ and
$\dg_t(x)=\dg_t(\g^{-1}_t(y))=\frac{d}{ds}(\g_{t+s}\circ\g^{-1}_t(y))=\Gamma
_R(y,
t)$.
We
thus have a parametrization of the canonical relation,
$$N^*\cM'=\Bigl\{(\g^{-1}_t(y),(D_x\g)^*(\eta);y,\eta): y\in\R^n, t\in\R,
\eta\perp\Gamma_R(y,t))\Bigr\},
\tag 6.5$$
which is favorable for analyzing the projection $\pr$. Indeed the
equivalence
of $(P)_R$ with the strong cusp condition follows immediately from the Lemma
in \S2.4.
%will follow from
%%part (a)
%of the following lemma.
%\footnote{ while
%implication
%$(P_k)_R\implies (C_k)_R$ will follow from part (b).
%}

\subheading{\bf 6.4. Pullback and commutator conditions}
The {\it bracket condition}  $(B)_R$ for families of curves (6.2)
states the linear independence of
vector fields $\hX{i}$, $i=1,\dots, n$ where
$\hX1=X_1$, $ \hX2=X_2$ and for $k=2,\dots, n$
$$\hX{k}:=X_k+\sum_{m=2}^{k-1}\sum_{I=(i_1,...,i_m)} a_{I,k}
[X_{i_1},[X_{i_2},...,[X_{i_{m-1}},X_{i_m}]\dots]]
\tag 6.6
$$
with universal coefficients $a_{I,k}$ which can be computed
 from the
coefficients of the Campbell-Hausdorff formula (\cite{40,Ch.V.5}, see also
the
exposition  in \cite{11}).
In particular
$$ \aligned
\hX1&=X_1, \quad\hX2=X_2,
\\
\hX3&=X_3-\frac16[X_1,X_2]
\\
\hX4&=X_4-\frac14[X_1,X_3]+\frac1{24}[X_1,[X_1,X_2]]
\\
\hX5&=X_5-\frac3{10}[X_1,X_4]-\frac1{10}[X_2,X_3]
+\frac1{15}[X_1,[X_1,X_3]]\\
&\quad\qquad+\frac1{30}[X_2,[X_1,X_2]]-\frac1{120}[X_1,[X_1,[X_1,X_2]]].
\endaligned
\tag 6.7
$$
See   \cite{56}, \cite{26} for the computation of the  vector
fields $\hX3$, $\hX4$ and their relevance for folds and cusps.

Assuming $(P)_R$ we shall now show that $(B)_R$ holds and how one can
determine the coefficients in (6.6).
By Taylor's theorem in the $s$ variable
$$\gamma_{s+t}\circ\gamma^{-1}_t
=\exp\Bigl(\phi(t,X_1,\dots,X_n,\dots)+s\psi(t,X_1,\dots,X_n,\dots)+O(s^2)\Bigr),
\tag 6.8
$$ and then, by an application of  the Campbell-Hausdorff formula
 (essentially  \cite{26, Eq. (6.4)}), we can
rewrite
this
as $$\exp\Bigl(O(s^2)\Bigr)\circ exp\Bigl(\phi + s\psi\Bigr).$$
From this it follows that
$$\Gamma_R(x,t)
%\frac{d}{ds}\bigl(\gamma_{s+t}^{-1}\circ\gamma_t(x)\bigr)
%\bigr|_{s=0}
=\psi(t,X_1,\dots,X_n,\dots)$$
and thus condition $(P)_R$ becomes the linear independence of
$\psi,\psi',...,\psi^{(n-1)}$. We will work modulo $O(s^2)+O(st^{n+1})$ and
so can
assume that there are only $n$ vector fields, $X_1,\dots,X_n$. Compute
$$\eqalign{\gamma_{s+t}\circ\gamma^{-1}_t
=&\exp(\sum_{i=1}^n(s+t)^iX_i)\circ \exp(-\sum_{i=1}^nt^iX_i)\cr
=&\exp((\sum t^iX_i +s\sum it^{i-1}X_i)+O(s^2))\circ \exp(-\sum t^iX_i)\cr
=&\exp(\sum_{i=1}^n(t+is)t^{i-1}X_i)\circ \exp(-\sum_{i=1}^n t^iX_i)\mod
O(s^2)\cr
=&\exp(B)\circ exp(A)}$$
with $A=-\sum_{i=1}^n t^iX_i$ and $B=\sum_{i=1}^n (t+is)t^{i-1}X_i$.
Now, the explicit  Campbell-Hausdorff formula (see \cite{40})  can be
written as
$$\align
\exp&(B)\circ \exp(A)
\\&= \exp
\Bigl(A+B+\frac12[A,B]+\sum_{m=3}^\infty
\sum_{\buildrel{I=(i_1,...,i_m)}\over{\in
\{1,2\}^m}}c_I \ad(C_{i_1})...\ad(C_{i_{m-1}})(C_{i_m})\Bigr)\cr
\\&=\exp
\Bigl(A+B+\frac12[A,B]+\sum_{m=3}^\infty\sum_{\buildrel{J=(j_1,...,j_{m-2})}
\over{\in
\{1,2\}^{m-2}}}\tilde c_J \ad(C_{j_1})...\ad(C_{j_{m-2}})([A,B])\Bigr)
\tag 6.9
\endalign
$$
where $C_1=A, C_2=B$ and $$\tilde c_J=c_{(J,1,2)}-c_{(J,2,1)}.
$$
The first few terms are
given by
$$\align 
A+B+&\frac12 [A,B]+\frac1{12}[A,[A,B]]-\frac1{12}[B,[A,B]]
\\&-\frac1{48}[A,[B,[A,B]]]-\frac1{48}[B,[A,[A,B]]]\dots.
\tag 6.10
\endalign$$
For notational convenience, we let the sum start at  $m=2$ instead of $m=3$
and set
$\tc_{\emptyset }=1/2$, and for the higher coefficients we get
$\tc_{(1)}=-\tc_{(2)}=1/12$ and $\tc_{(1,2)}=\tc_{(2,1)}=-1/48$. These are
enough to calculate
the coefficients in
$(B)_R$ in dimensions less than or equal to five
which is the situation  corresponding to at most
$S^+_{1,1,1,0}$ (strong swallowtail) singularities.

\par Returning to $(P)_R$, since we have $C_1=A, C_2=B$, we can use the
Kronecker delta notation  to write
$C_j=(-1)^{j}\sum_{i=1}^n(t+\delta_{j2}is)t^{i-1}X_i$.
%Denoting by ``$\sim$'' equality $\mod O(s^2)$ we obtain
%\pagebreak
Now
$$%\align
\gamma_{s+t}\circ\gamma^{-1}_t
=
\exp\Bigl(A+B+\sum_{m=2}^\infty\sum_{\buildrel{J=(j_1,...,j_{m-2})}
\over{\in
\{1,2\}^{m-2}}}\tilde c_J \ad(C_{j_1})...\ad(C_{j_{m-2}})([A,B])\Bigr)
+O(s^2)
$$
which modulo $O(s^2)$
 is equal to
$$\align&\exp\Bigl(-\sum_{i=1}^nt^iX_i+\sum_{i=1}^n(t^i+ist^{i-1})X_i\cr
&\qquad-\sum_{m=2}^\infty
\sum_{\buildrel{J=(j_1,\dots,j_{m-2})}\over{\in\{1,2\}^{m-2}}}
\tc_J
\ad\bigl((-1)^{j_1}\sum_{i_1}(t+\delta_{j_12}i_1s)t^{i_1-1}X_{i_1}\bigr)...
\cr
&\quad\dots\quad
\ad\bigl((-1)^{j_{m-2}}\sum_{i_{m-2}}(t+\delta_{j_{m-2}2}i_{m-2}s)
t^{i_{m-2}-1}X_{i_{m-2}}\bigr)\cr
&\qquad\cdot\bigl(\Bigl[\sum_{i_{m-1}=1}^n -t^{i_{m-1}-1}X_{i_{m-1}},
\sum_{i_m=1}^n(t+i_ms)t^{i_m-1}X_{i_m}\Bigr]\bigr)\quad\Bigr),
\endalign
$$
which, again modulo $O(s^2)$, is equal to
$$\align
&\exp\Bigl(\phi(t,X_1,\dots,X_n)+s\Bigl[\sum_{i=1}^nit^{i-1}X_i
-\sum_{m=2}^\infty
\bigl(\sum_{\buildrel{J}\over{\in\{1,2\}^{m-2}}}(-1)^{\sum_{l=1}^{m-2}j_l}
\tc_J\bigr)\times\cr
&\qquad\sum_{i_1,...,i_{m-2}}\quad\sum_{i_{m-1}<i_m}(i_m-i_{m-1})\cdot
\ad(X_{i_1})\cdot\dots\cr
&\qquad\quad\dots\cdot\ad(X_{i_{m-2}})\bigl([X_{i_{m-1}},X_{i_m}]
\bigr)t^{-1+\sum_{l=1}^mi_l}\Bigr]\Bigr).
\endalign
$$
From this we obtain
$$\align
\Gamma_R(x,t)=&\sum_{i=1}^nit^{i-1}X_i
-\sum_{m=2}^\infty
\bigl(\sum_{\buildrel{J}\over{\in\{1,2\}^{m-2}}}(-1)^{\sum_{l=1}^{m-2}j_l}
\tc_J\bigr)\times\cr
&\sum_{i_1,...,i_{m-2}}\,\sum_{i_{m-1}<i_m}(i_m-i_{m-1})\cdot
\ad(X_{i_1})\cdot\dots\cr
&\quad\dots\cdot \ad(X_{i_{m-2}})\bigl([X_{i_{m-1}},X_{i_m}]
\bigr)t^{-1+\sum_{l=1}^mi_l}\cr
:=&\sum_{i=1}^nit^{i-1}\hX{i}.
\endalign
$$
Since the $\tc_J$'s are known (\cf. \cite{40,Ch.V.5}, \cite{77})
this allows one to compute the $\hX{i}$'s and this shows that the condition
$(P)_R$
is equivalent with a bracket condition $(B)_R$ for some coefficients
$a_{I,k}$.

%\vfil\eject
To illustrate this, we restrict  to $n\le 5$ and  to get a manageable expression
we work 
$\mod O(t^5)$
and use (6.10);
 the expression for   $\Gamma_R(x,t)$  becomes then
$$\align
&\sum_{i}it^{i-1}X_i
-\frac12\sum_{i_1<i_2}(i_2-i_1)[X_{i_1},X_{i_2}]t^{i_1+i_2-1}\cr
&+\frac16\sum_{i_1}\sum_{i_2<i_3}(i_3-i_2)\cdot[X_{i_1},[X_{i_2},X_{i_3}]]
t^{i_1+i_2+i_3-1}\cr
&-\frac1{24}\sum_{i_1,i_2}\quad\sum_{i_3<i_4}(i_4-i_3)\cdot[X_{i_1},[X_{i_2},
[X_{i_3},X_{i_4}]]]t^{i_1+i_2+i_3+i_4-1}
\endalign
$$ which becomes
$$\align
&X_1+2tX_2+3t^2X_3+4t^3X_4+5t^4X_5\cr
&-\frac12[X_1,X_2]t^2-[X_1,X_3]t^3-\frac32[X_1,X_4]t^4-\frac12[X_2,X_3]t^4\cr
&+\frac16[X_1,[X_1,X_2]]t^3+\frac13[X_1,[X_1,X_3]]t^4
+\frac16[X_2,[X_1,X_2]]t^4\cr
&-\frac1{24}[X_1,[X_1,[X_1,X_2]]]t^4
\\
=&\hX{1}+2t\hX{2}+3t^2\hX{3}+4t^3\hX{4}+5t^4\hX{5}
\endalign
$$
where the $\hX{i}$ are given in (6.7).
Thus condition  $(B)_R$ in dimension $n\le 5$ is the linear independence of
the $\hX{i}$ for $0\le i\le n-1$.

\subheading{6.5. Curves on some nilpotent groups}

Let $G$ be an $n$ dimensional nilpotent Lie group with Lie algebra $\fg$.
 Let
$\g:\R\rta G$
be a smooth curve and define
$$G_R(t)=(DR_{\g(t)})^{-1}(\g'(t)),$$
where
$DR_g$ denotes the differential of right-translation by $g\in G$.
Note
that
 $G_R:\R\rta T_0 G=\frak \fg $ defines a curve in the Lie algebra $\fg$.

\proclaim{Lemma} The pullback condition $(P)_R$ for the family of curves
$t\mapsto  x\cdot \gamma(t)^{-1}$ is satisfied if and only if the vectors
$G_R(t), G_R'(t), \dots, G_R^{(n-1)}(t)$ are linearly independent
everywhere.
\endproclaim

To prove this, compute
$$
\align
\Gamma_R(x,t)=&\frac{d}{ds}\Bigl(\g_{s+t}\bigl(x\cdot\gamma(t)\bigr)
\Bigr)|_{s=0}\cr
=&\frac{d}{ds}\bigl(x\cdot\g(t)\cdot\g(s+t)^{-1}\bigr)|_{s=0}\cr
=&-x\cdot\g(t)\cdot\g^{-1}(t)\cdot\g'(t)\cdot\g^{-1}(t)\cr
=&-x\cdot\g'(t)\cdot\g^{-1}(t)\cr
=&-x\cdot \bigl(DR_{\g(t)}^{-1}(\g'(t)\bigr)=-x\cdot G_R(t),
\endalign$$
from which the equivalence is obvious.

The condition that $G_R,\dots G_R^{(n-1)}$ be linearly independent came up
in work of
Secco \cite{65}, who proved under this condition  the sharp 
$L^{3/2}\to L^2$ boundedness result
for the convolution operator
$$\cR f(x)=\int f(x\cdot\gamma(t)^{-1}) \chi (t) dt$$
on the Heisenberg group $\Bbb H$ (thus $n=3$). For the model family
of cubics
$\g(t)=(t,t^2,\a t^3)$, one easily computes that
$G_R(t)=(1,2t,(3\a+\frac16)t^2)$, so
that her condition is satisfied if and only if $\a\ne-\frac16$.

We further
illustrate the Lemma above by analyzing a two-parameter family of
quartics on a
four-dimensional, three-step nilpotent group, which we denote $\bbM$, due to
its
relation with the Mizohata operator.
The Lie algebra $\fm$ of $\bbM$ is spanned by $Y_j, 1\le j\le 4$, satisfying
$$[Y_1,Y_2]=Y_3,\quad [Y_1,Y_3]=Y_4,$$
with all other commutators equal zero. Thus, $Y_1$ and $Y_2$
satisfy the same commutator relations as
real and imaginary parts of the operator
$\frac{\partial}{\partial x} + i\frac{x^2}2 \frac{\partial}{\partial y}$,
{\it cf.}
\cite{45}.

The
group multiplication is given by
$$\multline
\bigl(x_1,x_2,x_3,x_4)\cdot (y_1,y_2,y_3,y_4)= (x_1+y_1,
x_2+y_2, x_3+y_3+\frac12(x_1y_2-x_2y_1),
\\
x_4+y_4+\frac12(x_1y_3-x_3y_1)+\frac1{12}(x_1-y_1)(x_1y_2-x_2y_1)\bigr).
\endmultline
$$

 For $\a,\bet\in\R$, we define  curves $\gamma(t)=(t,s^2,\a t^3,\bet t^4)$
and
ask for which values of the parameters  the vectors $G_R(t), \dots,
G_R'''(t)$ are
linearly independent.

We derive this in two different ways: first by the above  Lemma
 and then using the bracket condition.
 To form
$G_R$, we first calculate the derivative of $R_y(x)=x\cdot y$,
acting on a tangent vector $X=(X_1,X_2,X_3,X_4)\in\frak m=T_0\M$ :
$$DR_y(X)=(X_1,X_2,X_3+\frac12y_2X_1
-\frac12y_1X_2,X_4+\frac{6y_3-y_1y_2}{12}X_1
+\frac{y_1^2}{12}X_2-\frac{y_1}2X_3).$$
Computing the inverse of this and applying it for $y=\gamma(t)=(t,t^2,\a
t^3,\bet t^4)$, one calculates

$$\align
G_R(s)&=\Bigl(DR_{\gamma(t)}\Bigr)^{-1}(\dot\gamma (t))
=\Bigl(DR_{\gamma(t)}\Bigr)^{-1}(1,2t,3\a t^2,4\bet t^3)
\\
&=(1,2t,(\frac{6\a
+1}2) t^2,(\a+4\bet+\frac16) t^3).
\endalign
$$
Thus, $G_R^{(i)}$, $i=0,\dots,3$ are linearly independent if and only if
$\a+\frac16\ne0$ and
$\a+4\bet+\frac16\ne 0$.

Alternatively  we may quickly rederive this by using the bracket condition
$(B)_R$ for
$n=4$. We have
$$\g(x,t)=x\cdot (t,t^2,\a t^3,\bet t^4)^{-1}=\exp
\bigl(t(-Y_1)+t^2(-Y_2)
+ t^3(-\a Y_3)+t^4(-\bet Y_4)\bigr)(x),$$
where $Y_1,...,Y_4$ is the above basis for
$\frak m$, so we have the representation as in (1.1) with
$$X_1=-Y_1,\quad X_2=-Y_2,\quad X_3=-\a Y_3,\quad X_4=-\bet Y_4$$
and thus condition $(B)_R$ says that the vector fields
$$ -Y_1,\, -Y_2,\, -\a Y_3-\frac16[-Y_1,-Y_2],\,
-\bet Y_4-\frac14[-Y_1,-\a Y_3]+\frac1{24}[-Y_1,[-Y_1,-Y_2]]
$$
are linearly independent, which is equivalent with the linear independence
of the vector fields
$Y_1$, $ Y_2$, $(\a+\frac16)Y_3$ and
$(\frac{\a}4+\bet+\frac1{24})Y_4$.

\enddocument